\RequirePackage{fix-cm}
\documentclass[11pt,a4paper]{amsart}

\usepackage[a4paper,
inner=40mm,
outer=40mm,
top=40mm,
bottom=40mm,
marginparsep=5mm,
marginparwidth=40mm,
]{geometry}

\pdfoutput=1
\pdfpageattr{/Group <</S /Transparency /I true /CS /DeviceRGB>>}                  

\usepackage{graphicx}
\usepackage{mathptmx}      
\usepackage[normalem]{ulem}

\usepackage[american]{babel}
\usepackage[utf8x]{inputenc}
\usepackage[T1]{fontenc}

\usepackage{enumerate}
\usepackage{todonotes}
\usepackage{yfonts}
\usepackage{colortbl}
\definecolor{Gray}{gray}{0.9}
\definecolor{White}{gray}{1.0}
\newcolumntype{g}{>{\columncolor{Gray}}c}
\usepackage{amsmath, amsfonts}
\usepackage{mathtools}
\usepackage{dsfont}
\usepackage{csquotes}
\usepackage{booktabs}
\usepackage{graphicx}
\usepackage{color}
\usepackage[colorlinks=true, allcolors=blue]{hyperref}
\usepackage[ruled]{algorithm2e}

\usepackage{tikz}
\usepackage{xspace}
\usepackage{verbatim}
\usetikzlibrary{arrows, decorations.markings}
\usetikzlibrary{shapes.misc}
\usetikzlibrary{shapes.multipart}
\tikzset{cross/.style={cross out, draw=black, minimum size=2*(#1-\pgflinewidth), inner sep=0pt, outer sep=0pt},
	cross/.default={7pt}}
\renewcommand{\phi}{\varphi}

\newcommand{\R}{\mathds{R}}
\newcommand{\Rm}{\mathds{R}^m}

\newcommand{\suchthat}{\,\mid\,}
\newcommand{\define}{\coloneqq}
\newcommand{\card}[1]{\lvert{#1}\rvert}
\newcommand{\abs}[1]{\lvert{#1}\rvert}
\newcommand{\norm}[1]{\lVert{#1}\rVert}
\newcommand{\actuator}{a}
\newcommand{\frequency}{{\operatorname{freq}}}
\newcommand{\torque}{{\operatorname{torq}}}
\newcommand{\fuel}{{\operatorname{fuel}}}
\newcommand{\emission}{E}
\newcommand{\emconstraint}{\textfrak{e}}

\newcommand{\domain}{{U_{\rm ad}}} 
\newcommand{\codomain}{{Y_{\rm ad}}} 
\newcommand{\feasible}{{F_{\rm fsb}}} 
\newcommand{\grid}{G}
\newcommand{\graph}{{\mathcal G}}
\newcommand{\umeasurements}{\widehat{U}}
\newcommand{\dmeasurements}{\widehat{F}} 

\newcommand{\dcCont}{\gamma_{\operatorname{dc}}}

\newcommand{\dcDisc}{\operatorname{DC}}
\newcommand{\SOL}{\operatorname{SOL}}
\newcommand{\SOLprelim}{\widetilde{\operatorname{SOL}}}
\newcommand{\OP}{\operatorname{OP}}
\newcommand{\kOP}{{k\operatorname{OP}}}
\newcommand{\numberF}{\delta}
\newcommand{\boxweight}{w_{\operatorname{b}}}
\newcommand{\edgeweight}{w_{\operatorname{e}}}

\newcommand{\nedc}{\text{NEDC}\xspace}
\newcommand{\rand}{\text{RANDOM}\xspace}
\newcommand{\euro}{\text{EURO}\xspace}

\newcommand{\ecu}{\text{ECU}\xspace}
\newcommand{\Cameo}{\texttt{Cameo}\xspace}
\newcommand{\MATLAB}{\texttt{MATLAB}\xspace}
\newcommand{\LOLIMOT}{\texttt{LOLIMOT}\xspace}
\newcommand{\TopExpert}{\texttt{TopExpert}\xspace}
\newcommand{\mbminimize}{\texttt{mbminimize}\xspace}

\newcommand{\ASCMO}{\texttt{ASCMO}\xspace}

\newcommand{\RF}{\texttt{RF}\xspace}
\newcommand{\IF}{\texttt{IF}\xspace}
\newcommand{\MT}{\texttt{MT}\xspace}
\newcommand{\RP}{\texttt{RP}\xspace}
\newcommand{\TG}{\texttt{TG}\xspace}
\newcommand{\PM}{\texttt{PM}\xspace}
\newcommand{\PN}{\texttt{PN}\xspace}

\newcommand{\AF}{\texttt{AF}\xspace}
\newcommand{\PI}{\texttt{PI}\xspace}
\newcommand{\PT}{\texttt{PT}\xspace}
\newcommand{\NX}{\texttt{NO$_\texttt{x}$}\xspace}
\newcommand{\CO}{\texttt{CO}\xspace}
\newcommand{\CT}{\texttt{CO$_\texttt{2}$}\xspace}
\newcommand{\HC}{\texttt{HC}\xspace}

\DeclareMathOperator\vol{vol}
\DeclareMathOperator\minimize{minimize}

\theoremstyle{definition}
\newtheorem{definition}{Definition}[section]
\newtheorem{remark}{Remark}[section]

\newenvironment{myitemize}{%
	\begin{list}{$\circ$}%
		{\setlength{\topsep}{0.5ex}%
			\setlength{\partopsep}{0mm}%
			\setlength{\parskip}{0ex}%
			\setlength{\parsep}{0mm}%
			\setlength{\itemsep}{0ex}%
			\setlength{\labelwidth}{4mm}%
			\setlength{\leftmargin}{0ex}%
			\addtolength{\leftmargin}{\labelwidth}%
			\addtolength{\leftmargin}{\labelsep}%
			\setlength{\itemindent}{0mm}}}%
	{\end{list}}

\begin{document}
	
	\title{Semi-automatically optimized calibration of \\ internal combustion engines\thanks{Research by M. Joswig is partially supported by Einstein Stiftung Berlin and Deutsche Forschungsgemeinschaft (EXC 2046: \enquote{MATH$^+$}, SFB-TRR 109: \enquote{Discretization in Geometry and Dynamics}, SFB-TRR 195: \enquote{Symbolic Tools in Mathematics and their Application}, and GRK 2434: \enquote{Facets of Complexity}).\\
			The work of T. Burggraf, M. Pfetsch and S. Ulbrich has been supported by
			Deutsche Forschungsgemeinschaft within GSC 233 \enquote{Computational
				Engineering}. M. Pfetsch and S. Ulbrich have also been supported by
			the SFB 805 \enquote{Control of Uncertainty in Load-Carrying Structures in
				Mechanical Engineering} and S. Ulbrich by the GSC 1070 \enquote{Energy Science
				and Engineering}.
		}
	}
	
\title[Calibration of internal combustion engines]{Semi-automatically optimized calibration of \\ internal combustion engines}

\author[Burggraf et al.]{Timo Burggraf}
\address{Bertrandt Technikum GmbH, 71139 Ehningen, Germany}
\curraddr{}
\email{timoburggraf@gmx.de}
\thanks{The work of T. Burggraf, M. Pfetsch and S. Ulbrich has been supported by
Deutsche Forschungsgemeinschaft within GSC 233 \enquote{Computational
	Engineering}. M. Pfetsch and S. Ulbrich have also been supported by
the SFB 805 \enquote{Control of Uncertainty in Load-Carrying Structures in
	Mechanical Engineering} and S. Ulbrich by the GSC 1070 \enquote{Energy Science
	and Engineering}.}

\author[]{Michael Joswig}
\address{Einstein Professor of Discrete Mathematics/Geometry, Technical University of Berlin, 10623 Berlin, Germany}
\curraddr{}
\email{joswig@math.tu-berlin.de}
\thanks{Research by M. Joswig is partially supported by Einstein Stiftung Berlin and Deutsche Forschungsgemeinschaft (EXC 2046: \enquote{MATH$^+$}, SFB-TRR 109: \enquote{Discretization in Geometry and Dynamics}, SFB-TRR 195: \enquote{Symbolic Tools in Mathematics and their Application}, and GRK 2434: \enquote{Facets of Complexity}).}

\author[]{Marc E. Pfetsch}
\address{Research Group Optimization, Technische Universit\"at Darmstadt, 64293 Darmstadt, Germany}
\curraddr{}
\email{pfetsch@mathematik.tu-darmstadt.de}

\author[]{Manuel Radons}
\address{Workgroup Discrete Mathematics/Geometry, Technische Universit\"at Berlin, 10623 Berlin, Germany}
\curraddr{}
\email{radons@math.tu-berlin.de}

\author[]{Stefan Ulbrich}
\address{Research Group Optimization, Technische Universit\"at Darmstadt, 64293 Darmstadt, Germany}
\curraddr{}
\email{ulbrich@mathematik.tu-darmstadt.de}

	\maketitle
	
	\begin{abstract}
		Modern combustion engines incorporate a number of actuators and sensors
		that can be used to control and optimize the performance and
		emissions. We describe a semi-automatic method to simultaneously measure
		and calibrate the actuator settings and the resulting behavior of the
		engine. The method includes an adaptive process for refining the
		measurements, a data cleaning step, and an optimization procedure. The
		optimization works in a discretized space and incorporates the conditions
		to describe the dependence between the actuators and the engine behavior
		as well as emission bounds. We demonstrate our method on
		practical examples.
	\end{abstract}
	
	\section{Introduction}
	\label{sec:intro}
	Due to the wish to save fossil fuels, stringent maximal emission limits and challenging customer preferences, modern internal combustion engines (ICEs) become more and more complex.
	Indeed, research towards the improved construction of combustion engines is highly relevant for reaching the \CT emission targets set by the European Union, cf.~\cite{TERM2016}.
	
	This engine development results in an increasing number of \emph{actuators} and \emph{sensors}.
	There are currently in the order of ten different actuators and sensors each.
	Examples for actuators include the amount of injected fuel, exhaust recirculation control, air valve angle, etc.
	Sensors measure, e.g., the temperature, maximal point of the cylinder, torque, exhaust emission, etc.
	The actuators allow to produce a certain torque and revolution frequency, which describe the two main requirements on the engine in usage.
	However, several different settings of the actuators can result in the same
	torque/revolution frequency combination.
	Moreover, their dependence is involved and not known exactly a priori, i.e., it has to be measured and approximated.
	In this paper we deal with the \emph{optimal engine calibration problem}, i.e., to efficiently approximate this dependence by few measurements and to choose optimal actuators settings.
	
	The engine calibration problem consists of determining a so-called
	\emph{engine manifold}, which determines for each torque/revolution frequency
	combination a corresponding setting for actuators. This manifold is usually 
	discretized and the resulting \emph{solution map} is hard-coded into the engine control unit (\ecu).
	The settings on the engine manifold have to be chosen in such a way that they are consistent across various
	torque/frequency combinations, i.e., the engine manifold should be continuous. 
	Moreover, the resulting settings need to obey several
	restrictions in order to avoid damage of the engine as well as bounds on the emissions produced.
	The emissions are measured with respect to so-called \emph{driving cycles}.
	These are certain prescribed changes in the torque/frequency settings over time which are supposed to resemble usage in practice.
	These driving cycles are applied to the engine on a measurement bench and
	the resulting emissions have to be bounded.
	
	The calibration process described above involves two main steps.
	First, the dependence between the actuator settings and the output has to be measured.
	A naive enumeration of a uniform grid of possible actuator settings and
	interpolation would require exponential time in the number of actuators.
	Therefore, one needs to design a process to perform the measurements in relevant areas in order to speed up the measurement process.
	Additionally, the actuator settings are continuously changed without waiting until a steady-state has been reached.
	This allows to save time, but also introduces measurement
	errors like hysteresis effects.
	Second, based on the so-obtained information about the engine behavior, one should optimally choose the 
	resulting engine manifold and solution map.
	More precisely, one should select actuator settings for the final solution map that obey the above-mentioned 
	restrictions and yield an approximation to an optimal solution map which is sufficiently close in the sense 
	that all prescribed emission targets are met.
	
	In this article we propose a new way to solve the engine calibration problem,
	which consists of the following contributions:
	\begin{myitemize}
		\item \emph{Adaptive meshing:} The density of measurements is adapted
		within areas where the measured function is sensitive with respect to its
		inputs, while keeping the density of measurements coarse where it is
		not. This leads to a more accurate representation of the engine manifold than
		with a uniform grid approach with a fraction of the measurements.
		This first part is based on the well-established \LOLIMOT (local linear model tree)~\cite{Schueler2000} partitioning scheme of the space of actuator settings,
		which we extend by more involved measurement-routing and grid-refinement schemes. 
		
		\item \emph{Data cleaning:} Before optimization, the measured data are
		cleaned by filtering out redundancies and noise.
		
		\item \emph{Discretization and Optimization:} 
		We discretize the space of measurands in a fashion that fits the format of lookup tables as they are stored in the engine control unit. Using discrete optimization techniques, we select among the measured (and cleaned) data such actuator settings that minimize fuel consumption of the engine while its pollutant emissions conform to current regulations. The selected settings are drivable in the sense that the actuator's variation speed is bounded in order to prevent engine damage.  
	\end{myitemize}
	
	We would like to stress that the algorithm described in this article is a \emph{tool for}
	the engineer, not a replacement. While it runs automatically once its
	parameters have been set, the setting of these parameters, e.g., the
	determination of the subset of actuators which is to be varied in a given
	situation, requires extensive knowledge and experience. In this sense it is
	a semi-automatic process.
	\smallskip
	
	This article is structured as follows.
	In Section~\ref{sec:state-of-the-art} we briefly review relevant literature.
	A high-level description of the mathematical problem is given in Section~\ref{sec:math}.
	From this, we arrive at the corresponding steps of our process.
	Section~\ref{sec:calibration} provides the details of our method.
	In Section~\ref{sec:avl} we present a practical case-study, and Section~\ref{sec:experiments} contains the experimental results.
	We close with a few remarks.

	\section{State of the Art}\label{sec:state-of-the-art}
	
	In this section we describe the state of the art of calibration methods for the optimization of ICEs.
	Many commercial and research products exist for measurement and calibration.
	In the following we give a brief overview of the main approaches.
        Any calibration process comprises the following two major components:
	\begin{enumerate}[(A)]
		\item\label{enum:obtain} Obtain a good approximation of the engine behavior, mathematically described in terms of some function.
		\item\label{enum:lookup} Use whatever (approximate) knowledge of that function to produce a set of lookup tables for the \ecu.
	\end{enumerate} 
	The lookup tables may be optimized for various goals, e.g., dynamic driving behavior or maximal performance.
	In the present work we are set on optimizing the fuel consumption while conforming to a set of emission constraints. 
	
	\subsection{A Na\"{\i}ve Measurement-Based Approach}
	A basic idea is to measure all actuator settings on a sufficiently fine uniform grid.
	For modern engines this approach is infeasible for two reasons. 
	First, the size of the grid increases exponentially with the number of actuators.
	This would increase the number of measurements ---and thus also the measurement time--- beyond any feasible bounds.
	Second, even if such a comprehensive measurement were possible, it is computionally infeasible to optimize over such large input.
	So neither Step~\eqref{enum:obtain} nor Step~\eqref{enum:lookup} can be realized in this way.
	
	\subsection{The Model-Based Approach}
	For Step~\eqref{enum:obtain} in a model-based calibration the measurements are used to fit a given model to the engine behavior. 
	Once this fitting process is completed, in Step~\eqref{enum:lookup} lookup tables for the \ecu which are optimal with respect to the so-obtained model and a given set of objectives are computed via standard techniques from nonlinear optimization such as steepest descent methods, cf.~\cite[p. 542]{Isermann}.
	A software package which does so in an automated fashion is the \emph{Model-Based Calibration Toolbox} for \MATLAB \cite{MatlabECU2018}.
	Leading model-based approaches employ physical models, the training of neural networks via the measurements, and statistical machine learning. 
	An implementation that combines the first two kinds is, e.g., the software package \mbminimize~\cite{Knoedler2003}.
	Statistical machine learning is used, e.g., in \ASCMO~\cite{KRUSE2010739}.
	
	Physical models are usually described in terms of smooth functions. 
	Similarly, most neural networks also model smooth functions.
	For instance, perceptrons are commonly composed of sigmoid functions and thus smooth, cf.~\cite[pp.~103--111]{Isermann2003}.
	This restriction to mathematical modeling via smooth functions inherently imposes severe limitations to any model-based calibration of ICEs.
	The reason is that the functions necessary to describe an ICE well enough exhibit strong nonlinearities and may even be nondifferentiable.
	In the former case, the approximation error can be bounded globally, but locally the 
	approximating functions are likely to exhibit strong oscillations that do not describe the 
	underlying engine behavior very well; cf.~\cite[p.~105]{Isermann2003}. 
	At a nondifferentiability of a continuous function the role of the local linear approximation, which is the derivative, is taken over by local piecewise linear approximations~\cite[p.~67ff]{scholtes2012introduction}.
        In this case any sufficiently good approximation usually requires an exponential number of local models, one for each linear piece.
	As a consequence, it is common that either Step~\eqref{enum:obtain} does not describe the ICE well enough or both Step~\eqref{enum:obtain} and Step~\eqref{enum:lookup} are infeasible due to a combinatorial explosion, just like in the na\"{\i}ve approach.
	
	\subsection{LOLIMOT}
	Local linear model trees (\LOLIMOT) are arguably the most performant model-based calibration method to date; cf.~\cite[p.~93]{Isermann}, as well as~\cite{Martini2003}. 
	\LOLIMOT partitions the space of actuator settings as a dissection into cubical cells.
	In each cubical cell, the engine behavior is modeled 
	by a Gaussian function on the basis of a \enquote{central composite measurement point pattern}; cf.~\cite[Fig.~3.4.12]{Isermann}. 
	These local models are then stitched together into a global one.
	If a local model fails to produce a sufficiently close approximation, e.g., due to a
	strongly nonlinear behavior of the engine, the corresponding cubical cell is split along one 
	of its axes and the engine behavior is modeled via a central composite design on each sub-cell.
	This process is iterated until a sufficient model quality is reached globally.
	
	The key contribution of \LOLIMOT is to introduce adaptive meshing to the engine calibration process. 
	One drawback of the method is that only one output value is modeled at a time. 
	The engine model is then fused together from the component models for the individual measurands; cf.~\cite[p.~93]{Isermann}. 
	This approach may cause both holes in the modeled behavior as well as local redundancies 
	of data. 
	Further, the model produced by \LOLIMOT is a special type of radial basis function network; cf.~\cite[p.~143]{Isermann2010}.
	These types of neural nets require a locally homogeneous covering of the actuator space by measurements; cf.~\cite[p.~111]{Isermann2003}.
	This then necessitates the aforementioned fixed measurement patterns within each cell.
	But any such fixed pattern is again, like in the previous approaches, subject to a combinatorial explosion in high dimensions.
	The adaptive meshing with fixed local pattern is better than the uniform grid used in the na\"{\i}ve approach, but only by a constant factor.
	This is the second drawback.
	
	Our new method for engine calibration refines \LOLIMOT by simultaneously considering all measurands.
	This employs a more involved grid refinement scheme and randomized measurement routing; cf.\ Section~\ref{sec:iteration-step}.
	In this way we can overcome conceptual limitations of \LOLIMOT.
	
	\subsection{Data Boundaries}\label{sec:data-boundaries}
	\label{subsec:boundaries}
	Throughout its operation the measurand values of an ICE have to stay within certain boundaries.
	Some of these boundaries are set in place to avoid destruction, e.g., for cylinder pressures or critical device temperatures.
	Others are induced, e.g., by emission and noise regulations. 
	In practice they are obtained automatically via established software packages such as \Cameo~\cite{Gschweitl2001}, \TopExpert~\cite{TOPexpert2004}, or the \MATLAB toolbox \LOLIMOT (local linear model tree)~\cite{Schueler2000}.
	Throughout we will assume that all such boundary informations are already given.

	\section{Mathematical Problem Description}
	\label{sec:math}
	
	\noindent
	The calibration of an internal combustion engine is the procedure to derive
	an engine manifold, which is optimal with respect to some predefined objective. 
	To reduce the complexity of the model and to
	enhance practical implementation, these manifolds are discretized to obtain
	solution maps (see Section~\ref{sec:discretization} below).
	The calibration problem posed informally in~\cite[p.~250]{MST2018} asks for an engine manifold 
	that minimizes fuel consumption, while conforming to a number of emission constraints.
	In this section we will present two mathematical formalizations of the latter, one
	idealized continuous version and its discretization whose output fits the format of 
	the lookup tables for the ECU.  
	
	In our setting, knowledge of the engine behavior with respect to variations of its actuator settings is obtained by means of physical experiments on a test bench. In addition to the revolution frequency, typical actuators include the injected fuel quantity, the injection angle, or the valve pressure.
	The generated torque of the engine is a measurand. 
	
	Technically there is a relevant distinction between direct and controlled actuators. Direct actuators 
	such as injected fuel quantity, injection angle, or valve pressure can be set directly on the engine, while controlled actuators are set indirectly. For example, the revolution frequency is a controlled actuator that is regulated via a brake on the engine shaft. However, in our mathematical model this distinction is not relevant.
	
	The aforementioned side constraints include limits on pollutant emission as well as physical requirements such as engine temperature limits.
	They necessitate that not only torque, but several other output values of the engine are measured as well. 
	A realistic engine model features $m\geq 8$ actuators and $n\geq 14$ measurands.
	In our mathematical model we represent the relation between the setting of
	$m$ actuators and $n$ sensor values by a function
	\[
	F\colon \R^m \to \R^n\,.
	\]
	Throughout we make the fundamental assumption that $F$ is continuous, but not necessarily differentiable everywhere.
	Further, we will assume the actuator values to be restricted to a box $\domain\subset \R^m$, which we call the \emph{admissible domain}. 
	It was already noted in Section~\ref{sec:data-boundaries} that throughout this work
	we assume $\domain$ to be already given.
	The noncritical sensor values define another box $\codomain\subset \R^n$, the \emph{admissible range}; cf.\ Section~\ref{subsec:boundaries}.
	The exact definitions of $\domain$ and $\codomain$ will be stated in the subsequent Section~\ref{sec:init}.
	The \emph{feasible space} of $F$ is the set
	\[
	\feasible \define \{ (u,y)\in \domain\times\codomain \suchthat y=F(u) \} \ \subseteq \ \R^m\times\R^n\,.
	\]
	The actuators correspond to the coordinates of $u$; examples are the revolution frequency and the amount of fuel injected.
	Typical sensor values, i.e., coordinates of $F(u)$, include the torque and the emission of carbon monoxide.
	
	We denote by \enquote{$\frequency$} the index of the actuator for revolution frequency and by \enquote{$\torque$} the index of the torque sensor values.
	Then we call 
	\[
	\OP \define \{(u_\frequency,y_\torque) \suchthat (u,y)\in\feasible\}
	\]
	the \emph{operation field} of $F$.
	The operation field represents the behavior of the engine with respect to revolution frequency and torque.
	While our methods are more general, we focus on this particular pair of actuator and sensor values in our analysis.

	\subsection{Continuous Optimization Problem}\label{sec:continuous-opt}
	
	It is important to understand that the actual optimization problem is inherently discrete, since the desired output is a solution map for the ECU.
        These lookup tables of actuator values for given frequency and torque combinations have a given finite length.
        For the sake of a concise exposition, however, we now describe an idealized continuous optimization problem which has the actual optimization
	problem that we want to solve as a natural discretization.
        While similar optimization problems must be behind all known approaches to ICE calibration, we are not aware of any complete description in the available literature.
        Our model below is intended to fill this gap.
        
	\subsubsection*{Optimization Space and Drivability}
	
	The feasible region of our continuous optimization problem is given by all engine manifolds.
	Each such manifold is given as the image of a map $M \colon \OP \to \R^m$ that assigns actuator settings to a given frequency and torque value pair in the operation
	field, i.e.,
	\begin{align}\label{eq:M-constraint}
        [F(M(u_\frequency, y_\torque))]_\torque \ = \ y_\torque\quad\text{for all } (u_\frequency,
	y_\torque)\in \OP \enspace.
	\end{align}
	A basic requirement is that these maps are continuous.
	
	Moreover, there are vital additional conditions to consider.
	Any solution to the engine calibration problem must be drivable in the sense of~\cite[p.~258]{MST2018}:
	Varying actuators too fast might damage the engine.
	Therefore in the (continuous) final solution map the variation speed of every actuator is bounded by constants $\Delta_\actuator$ for $\actuator \in \{1, \dots,
	m\}$.
        For a given map $M$ and actuator~$\actuator$, the corresponding \emph{drivability constraint} 
        is that for all $(u_\frequency, y_\torque)$, $(u_\frequency',
        y_\torque') \in \OP$ with $u \define M(u_\frequency,y_\torque)$ and
        $u' \define M(u_\frequency',y_\torque')$ the following has to hold
	\begin{align}\label{eq:driv-constraint}  
	\abs{u_a-u_a'} \, = \,
	\abs{[M(u_\frequency,y_\torque)]_a-[M(u_\frequency',y_\torque')]_a}
	\,  \leq \, \Delta_\actuator \cdot \norm{(u_\frequency, y_\torque)-(u_\frequency', y_\torque')}\,,
	\end{align}
	where $\norm{\cdot}$ is some norm.
	As a consequence, the map $M_\actuator:\OP\to\R$ is Lipschitz continuous with constant $\Delta_\actuator$.
	We define $\Omega$ as the set of all maps~$M$ which are feasible in the sense that they obey the drivability constraint with respect to each
	actuator, more precisely,
	\[
         \Omega \ \define \ \bigl\{ M \colon \OP \to \R^m \colon\, M \text{ satisfies
        \eqref{eq:M-constraint} and moreover \eqref{eq:driv-constraint}
	for all $\actuator \in \{1, \dots, m\}$}\bigr\}.
	\]
	
	\subsubsection*{Driving Cycle and Emission Constraints}
	
\begin{figure}[t]
		\includegraphics[width=0.48\textwidth]{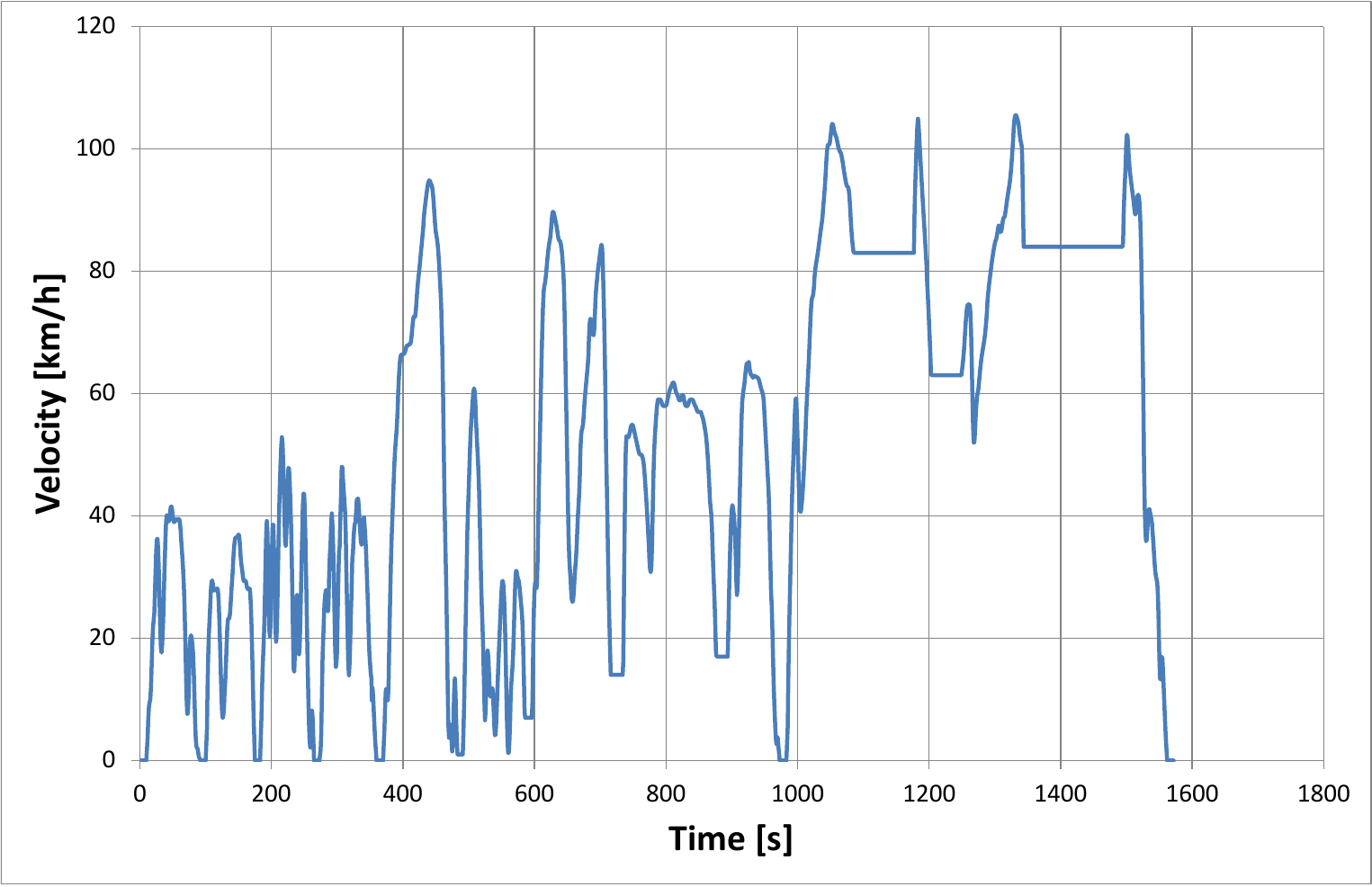}
		\hfill
		\includegraphics[width=0.48\textwidth]{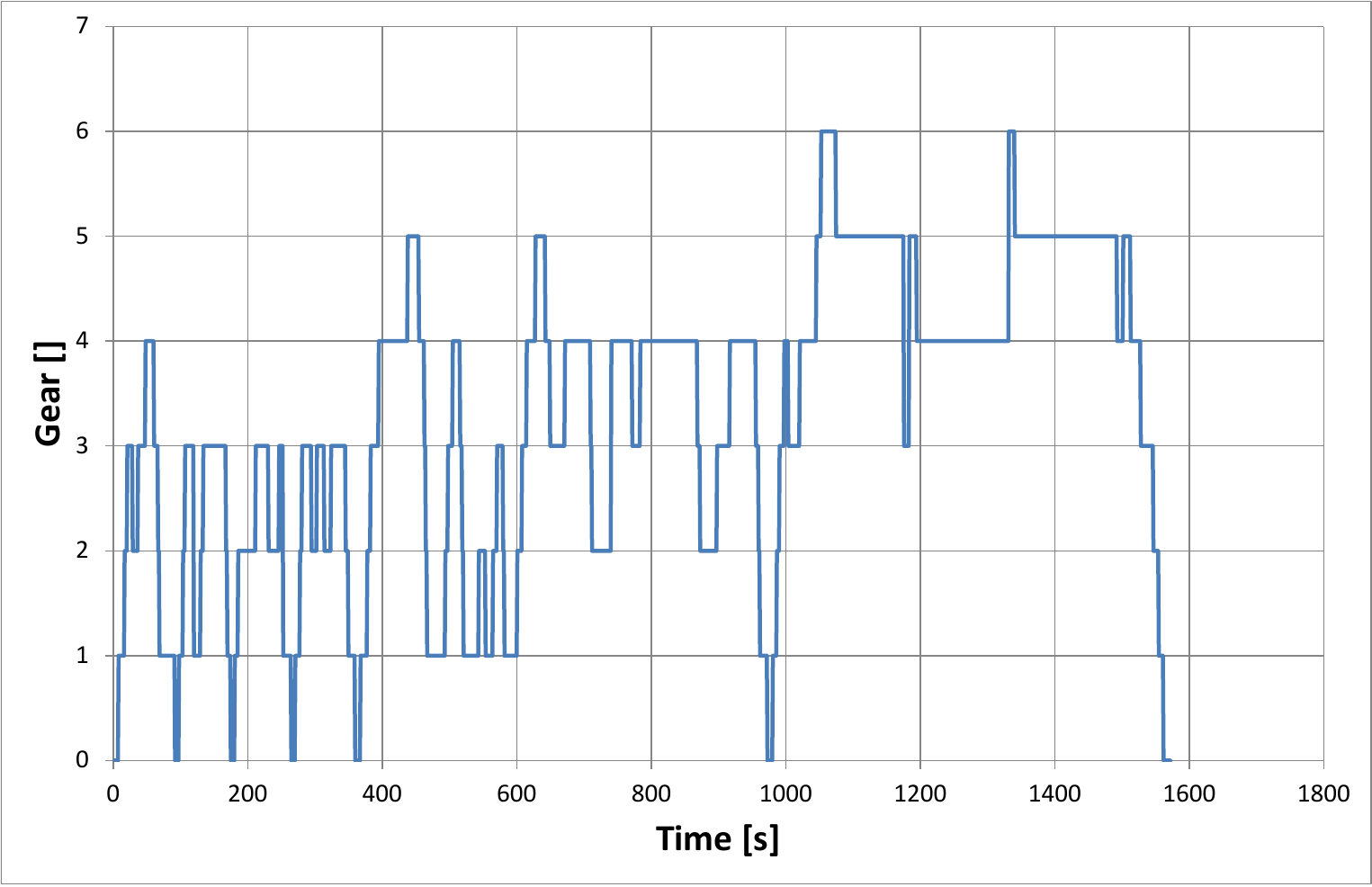}
		\caption{Operational profile of \rand driving cycle.}
		\label{fig:rnd-profile}
	\end{figure}

	In accordance with current government regulations and common test cycles, the engine behavior is optimized with respect to pre-defined scenarios, known as driving cycles.
	In our continuous model a \emph{driving cycle} is a time-parametrized curve
	\[
	\dcCont \colon [0,1] \to \OP\,,
	\]
	whose purpose is to simulate the phases of acceleration and constant speed of real-world driving patterns.
	Driving cycles are given indirectly in the form of operational profiles which map a time to a gear/velocity combination; cf.\ Figure~\ref{fig:rnd-profile}. 
	Taking into account the weight of the car, its drag coefficient,
	specific tire friction, etc., every such gear/velocity combination can be mapped 
	to a revolution frequency/torque combination in the operation field. 
	
	For the sake of the simplicity of exposition, here we will assume that the calibration is
	performed with respect to a fixed driving cycle.
	Yet it is also natural to take the combination of several driving cycles into account, and our approach also covers this slightly more general situation.
	
	An important constraint prescribed by regulations is to bound the resulting emissions along the driving cycle.
	Emission pollutants include carbon monoxide (\CO), hydrocarbons (\HC), nitrogen oxides (\NX) as well as
	particulate matter (\PM) and number (\PN). We denote by $\emission$ the index
	set that corresponds to pollutant emissions and by $\emconstraint_p$ the
	emission limit for pollutant~$p$ over the driving cycle.
	
	Consider an engine manifold map $M \in \Omega$ and a pollutant~$p \in \emission$.
	The frequency and torque values along $\dcCont$ produce actuator settings $M(\dcCont(t))$, which result in output values $F(M(\dcCont(t)))$, including the pollutant~$p$.
	The integral of these values must satisfy the \emph{emission constraint} for all $p \in \emission$:
	\begin{align}\label{eq:em-constraint}
	\int_{0}^{1} [F(M(\dcCont(t)))]_p \,\cdot\, \lVert \dot{\gamma}_{\operatorname{dc}}(t)\rVert\,\text{d}t\ \leq\ \emconstraint_p\,.
	\end{align} 
	The factor $\lVert \dot{\gamma}_{\operatorname{dc}}(t)\rVert$ accounts for the fact that $\dcCont$ is not necessarily normalized, i.e., the speed of acceleration, deceleration, etc., varies.
	Table~\ref{tab:euro-norms} shows the diesel engine emission constraints for \euro norms 3-6c (E3-E6c).
	
	\subsubsection*{The optimization formulation}
	
	Let \enquote{$\fuel$} be the index of the actuator for the injected fuel quantity. 
	It is possible that one wants to optimize fuel consumption over a different curve than the driving cycle $\dcCont$.
	We thus define a second curve in the operation field
	\[
	\phi \colon [0,1] \to \OP\,.
	\]
	Using the terminology developed so far, the continuous optimization problem is
	\begin{align}\label{eq:continuous-opt}
	\minimize_{M \in \Omega}\ &  \int_0^1 \, [M(\phi(t))]_\fuel \,\cdot\, \lVert \dot{\phi}(t)\rVert\,\text{dt}
	\\
	\text{subject to}\  & \int_{0}^{1} [F(M(\dcCont(t)))]_p \,\cdot\, \lVert \dot{\gamma}_{\operatorname{dc}}(t)\rVert\,\text{d}t\ \leq\ \emconstraint_p \quad\text{for all } p \in \emission\,.\notag
	\end{align}
	That is, the optimization goal is to minimize the total consumed fuel with respect to the
	time parametrized curve $\phi$ in $\OP$, while for all pollutants $p$ the integral emission along~$\dcCont$ is bounded by the prescribed constant $\emconstraint_p$. 
	Implicitly, problem \eqref{eq:continuous-opt} is subject to the
        driviability constraint, due to the fact that all manifold maps~$M \in \Omega$ must conform to inequality \eqref{eq:driv-constraint}. 
	In the discrete formulation of the engine calibration problem this dependence will be made explicit.
	We further remark that the curves $\phi$ and $\dcCont$ may coincide, but they do not have to. 
	More elaborate curves $\phi$ may be useful, e.g., for controlling the fuel consumption also outside the driving cycle, while ignoring the fringes of the operation field. 
	
	\renewcommand{\arraystretch}{1.3}
	\begin{table}	
		\caption{\label{tab:euro-norms} Selection of \euro emission constraints for passenger cars with compression ignition engine.}
		\small
		\centering
		\begin{tabular*}{\textwidth}{@{}l@{\extracolsep{\fill}}lrrrrr@{}}\toprule
			Em. & Unit &  E3 (2001) &  E4 (2006) &  E5a (2011) & E6b (2015) & E6c (2018) \\ \midrule
			\CO & $\text{mg}\,/\,\text{km}$ & 2300 & 1000 & 1000 & 1000 & 1000 \\
			\HC & $\text{mg}\,/\,\text{km}$ & 200 & 100 & 100 & 100 & 100 \\
			\NX & $\text{mg}\,/\,\text{km}$ & 150 & 80 & 60 & 60 & 60 \\
			\PM & $\text{mg}\,/\,\text{km}$ & -- & -- & 5 & 4.5 & 4.5 \\
			\PN & $$1\,/\,\text{km}$$ & -- & -- & -- & $6\cdot 10^{12}$ & $6\cdot 10^{11}$ 	\\\bottomrule
		\end{tabular*}
	\end{table}

	\subsection{Discrete Optimization Problem} \label{sec:discretization}
	
	To obtain a finite-dimensional problem, the model is discretized to yield
	\emph{(characteristic) engine maps}. We consider combinations of~$k$
	revolution frequencies and $k$ torque demands and subdivide the operation
	field $\OP$ of~$F$ into $k^2$ congruent rectangles $\OP_{ft}$,
        where $f$ and $t \in[k] \define \{1, \dots, k\}$ denote the rectangle's frequency and torque coordinate,
	respectively.
        This corresponds to the technical requirement that engine maps are given as $k {\times} k$-matrices, which
	for each combination yield the corresponding value of a particular
	actuator.  They are stored permanently in the engine control unit.
	A \emph{solution map} consists of a complete set of engine maps, one for each actuator.
	In this way a solution map yields a discretization of the continuous solution map
	$M$ described in Section~\ref{sec:continuous-opt}.
	
	Then, in our terminology, the solution map takes as input a frequency and
	torque pair $(f,t)\in [k]\times[k]$ and yields as output an admissible
	actuator setting.
	The latter is a point $u\in\domain$ with $(u,F(u))\in\feasible$
	such that the pair $(u_\frequency, y_\torque)$ lies in the rectangle of the
	discretized operation field corresponding to the input coordinates $(f,t)$.
	
	Our goal is to obtain a solution map which is optimal with respect to a
	given objective, while conforming to several constraints.  Since typical
	constraints are continuous but nonlinear and the solution map itself is
	discrete, this optimization problem belongs to the wide class of
	mixed-integer/discrete nonlinear optimization problems, which are often
	difficult to handle.
	
	Splitting the engine calibration into Steps~\eqref{enum:obtain} 
	(data acquisition) and~\eqref{enum:lookup} (computation of a solution map), cf.
	Section~\ref{sec:state-of-the-art}, our goal can be formulated as follows.
	In Step~\eqref{enum:obtain} obtain, via actual measurements of the engine, a finite set 
	\[
          \dmeasurements \define \{(u^1,y^1),(u^2,y^2),\dots,(u^\numberF,y^\numberF)\}
	\]
	of $\numberF$ points in $\feasible$, i.e., a set of actuator settings $u^q\in\domain$ such that $y^q=F(u^q)\in\codomain$, which is a sufficiently good representation of $F$.
        Here \enquote{sufficiently good} means that we can, in Step~\eqref{enum:lookup}, extract from $\dmeasurements$ a solution
	map that conforms to discretized versions of the drivability and emission constraints 
	presented in Section~\ref{sec:continuous-opt}. 
	In this sense, a solution map is a sufficiently good discrete approximation to an optimal
	engine map $M$ 
	of the continuous optimization problem~\eqref{eq:continuous-opt}.
	
	Since $\dmeasurements$ is obtained via actual measurements, it is called the \emph{data set}.
	The elements 
	\[
          d^q\define(u^q,y^q)
        \]
	of $\dmeasurements$ are called \emph{data points}.
	The final selection of the data points is the result of a cycle of measurements and optimization steps, which are discussed in the sections below.
	We will assume that none of the points in the data set lies on the boundary of any rectangle $\OP_{ft}$.
	In this case, the discretization of the operation field into $k^2$ rectangles partitions $\dmeasurements$ into sets 
	\[
	S_{ft}\define \bigl\{ (u,y) \in \dmeasurements \suchthat (u_{\frequency},y_{\torque}) \in\OP_{ft} \bigr\}\,.
	\]	
	We call each set $S_{ft}$ a \emph{stack}, and the entire partition 
	\[
	\kOP \define \{S_{ft} \suchthat f,t\in[k]\}
	\]
	is the \emph{$k$-operation field} of $F$ with respect to $\dmeasurements$.
	In the solution map each entire rectangle $\OP_{ft}$ will be represented by a single measurement $d_{ft}\in S_{ft}$.
	Hence, the discrete analogue to $\Omega$, the set of all manifold maps, is  
	the set $\Omega_k$ of all maps $M_k \colon [k]\times [k] \to \R^m\times \R^n$, such
        that $M_k(f,t)\in S_{ft}$, i.e., select exactly one data point from each
        stack, and satisfy a discrete analogue of the drivability constraint
        \eqref{eq:driv-constraint}, see \eqref{eq:driv-constraint-disc}
        below. Thus,
	\[
	\Omega_k \ \define \ \bigl\{ M_k \colon [k]\times [k] \to \R^m\times \R^n \colon\,
        \text{$M_k(f,t)\in S_{ft}$ for all $f,t\in[k]$ and $M_k$ satisfies
	\eqref{eq:driv-constraint-disc}}\bigr\} \enspace.
	\] 
	The elements of $\Omega_k$ are called solution maps. 
	Each solution map $M_k\in \Omega_k$ has precisely $k^2$ values $M_k(f,t)\in
	S_{ft}$, one per stack. A solution map is uniquely determined by the
	data points $d_{ft}=M_k(f,t)$, $f,t\in[k]$, and we call $d_{ft}$
        the \emph{representatives} of the $k^2$ stacks $S_{ft}$
	corresponding to $M_k$.
	The final calibration solution is a solution map $\SOL\in\Omega_k$ which is
	optimal with respect to the given objectives, i.e., minimization of
        fuel consumption, subject to emission regulations and drivability. 
	The final solution will be picked by solving the optimization problem \eqref{eq:discrete-opt} below.
	
	\subsubsection*{Discrete drivability constraint}
	
	Due to the uniform discretization of the $k$-operation field, \emph{a discrete drivability constraint} merely has to bound the difference between actuator settings of representatives of neighboring rectangles $\OP_{ft}$.
        That is, with $d_{ft}=M_k(f,t)$, the discrete analogue of \eqref{eq:driv-constraint} is
	\begin{align}\label{eq:driv-constraint-disc}
	\abs{[d_{ft}]_\actuator - [d_{gs}]_\actuator} \ = \ 
        \abs{[M_k(f,t)]_\actuator-[M_k(f,t)]_\actuator}\ \leq\ \Delta_\actuator \qquad \text{for all } \actuator\in[m]\,,
	\end{align}
	and all tuples $(f,t), (g,s)\in [k]\times [k]$, where either $g = f \pm 1$ or $s = t\pm 1$. 
	
	\subsubsection*{Discrete emission constraint}

        To discretize the emission constraint \eqref{eq:em-constraint}, we
	replace an engine map $M\in\Omega$ by its discrete counterpart
        $M_k\in\Omega_k$ and obtain instead of the curve integral along the
	curve $\dcCont$ a weigthed sum. In fact,
	to each rectangle $\OP_{ft}$ we associate a weight $\omega_{ft}$.
	This weight is set to zero if the intersection of $\OP_{ft}$ with the image of the curve $\dcCont$ is empty.
	Otherwise $\omega_{ft}$ is a positive value that reflects the resistance time, i.e., the duration of the curve
        $\dcCont$ staying in the rectangle $\OP_{ft}$.
	The function
	\begin{equation}\label{eq:dcDisc}
          \dcDisc \colon [k]\times[k]\to\R_{\geq 0},\ (f,t)\mapsto \omega_{ft}
	\end{equation}
	serves as a discrete analogue of the continuous driving cycle $\dcCont$.
	Note that the map $\dcDisc$ only records which parts of the operation field are met by  $\dcCont$
	for which duration, but it ignores the order in which this happens.
	
	The \emph{discrete emission constraint} is now given as
	\begin{equation}\label{eq:em-constraint-disc}
	\sum_{(f,t)\in [k]\times [k]} \omega_{ft}\cdot [M_k(f,t)]_p \ \leq\ \emconstraint_p \qquad\text{for all } p\in\emission\,.
	\end{equation}
	The practical driving cycles, for instance the New European Driving Cycle (\nedc) or Real World Driving Cycle (\rand), are given by operational profiles which map a time to a gear-velocity combination; cf.\ Figure~\ref{fig:rnd-profile} and Section~\ref{sec:continuous-opt}.  
	The weights $\omega_{ft}$ can be derived from these profiles. The revolution frequency at a certain time can be calculated directly from the current gear/velocity combination. 
	
	For the requested engine torque not only the speed but also the acceleration has to be taken into account, along with the car's mass, roll drag and air flow resistance. 
	For a schematic of a discretized driving cycle, cf.\ Figure~\ref{fig:frog}.
	
	\begin{figure}
		\centering
		
		\definecolor{gray1}{gray}{0.87}
		\definecolor{gray2}{gray}{0.6}
		\definecolor{gray3}{gray}{0.75}
		\begin{tikzpicture}[scale=1]
		
		
		\draw [<->,thick] (0,6.1) node (yaxis) [above] {$\torque$}
		|- (7.5,0) node (xaxis) [right] {$\frequency$};
		
		\draw (1,0) -- (1,5.6);
		\draw (2,0) -- (2,5.6);
		\draw (3,0) -- (3,5.6);
		\draw (4,0) -- (4,5.6);
		\draw (5,0) -- (5,5.6);
		\draw (6,0) -- (6,5.6);
		\draw (7,0) -- (7,5.6);
		
		\draw (0.5,0) -- (0.5,5.6);
		\draw (1.5,0) -- (1.5,5.6);
		\draw (2.5,0) -- (2.5,5.6);
		\draw (3.5,0) -- (3.5,5.6);
		\draw (4.5,0) -- (4.5,5.6);
		\draw (5.5,0) -- (5.5,5.6);
		\draw (6.5,0) -- (6.5,5.6);
		
		\draw (0,0.8) -- (7,0.8);
		\draw (0,1.6) -- (7,1.6);
		\draw (0,2.4) -- (7,2.4);
		\draw (0,3.2) -- (7,3.2);
		\draw (0,4) -- (7,4);
		\draw (0,4.8) -- (7,4.8);
		\draw (0,5.6) -- (7,5.6);
		\draw (0,0.4) -- (7,0.4);	
		\draw (0,1.2) -- (7,1.2);
		\draw (0,2) -- (7,2);
		\draw (0,2.8) -- (7,2.8);
		\draw (0,3.6) -- (7,3.6);
		\draw (0,4.4) -- (7,4.4);
		\draw (0,5.2) -- (7,5.2);
		
		\draw [fill=gray1] (0.5,0.8) rectangle (1,1.2);
		\draw [fill=gray1] (1,0.8) rectangle (1.5,1.2);
		\draw [fill=gray1] (1.5,0.8) rectangle (2,1.2);
		\draw [fill=gray3] (2,0.8) rectangle (2.5,1.2);
		\draw [fill=gray2] (2.5,0.8) rectangle (3,1.2);
		\draw [fill=gray2] (3,0.8) rectangle (3.5,1.2);
		\draw [fill=gray2] (3.5,0.8) rectangle (4,1.2);
		\draw [fill=gray3] (4,0.8) rectangle (4.5,1.2);
		\draw [fill=gray2] (4.5,0.8) rectangle (5,1.2);
		\draw [fill=gray3] (5,0.8) rectangle (5.5,1.2);
		\draw [fill=gray1] (5.5,0.8) rectangle (6,1.2);
		\draw [fill=gray1] (6,0.8) rectangle (6.5,1.2);
		
		\draw [fill=gray1] (1,1.2) rectangle (1.5,1.6);
		\draw [fill=gray1] (3,1.2) rectangle (3.5,1.6);
		\draw [fill=gray2] (4,1.2) rectangle (4.5,1.6);
		
		\draw [fill=gray1] (1,1.6) rectangle (1.5,2);
		\draw [fill=gray1] (3,1.6) rectangle (3.5,2);
		\draw [fill=gray1] (4,1.6) rectangle (4.5,2);	
		
		\draw [fill=gray1] (1,2) rectangle (1.5,2.4);
		\draw [fill=gray1] (1.5,2) rectangle (2,2.4);
		\draw [fill=gray1] (2,2) rectangle (2.5,2.4);
		\draw [fill=gray1] (2.5,2) rectangle (3,2.4);
		\draw [fill=gray1] (3,2) rectangle (3.5,2.4);
		\draw [fill=gray1] (3.5,2) rectangle (4,2.4);
		\draw [fill=gray1] (4,2) rectangle (4.5,2.4);
		
		\draw [fill=gray1] (1,2.4) rectangle (1.5,2.8);
		\draw [fill=gray1] (2,2.4) rectangle (2.5,2.8);
		\draw [fill=gray1] (3,2.4) rectangle (3.5,2.8);
		\draw [fill=gray1] (4.5,2.4) rectangle (5,2.8);
		
		\draw [fill=gray1] (1,2.8) rectangle (1.5,3.2);
		\draw [fill=gray1] (2,2.8) rectangle (2.5,3.2);
		\draw [fill=gray1] (3,2.8) rectangle (3.5,3.2);
		\draw [fill=gray1] (4.5,2.8) rectangle (5,3.2);
		
		\draw [fill=gray1] (1.5,3.2) rectangle (2,3.6);
		\draw [fill=gray1] (2,3.2) rectangle (2.5,3.6);
		\draw [fill=gray3] (2.5,3.2) rectangle (3,3.6);
		\draw [fill=gray1] (3,3.2) rectangle (3.5,3.6);
		\draw [fill=gray3] (3.5,3.2) rectangle (4,3.6);
		\draw [fill=gray3] (4,3.2) rectangle (4.5,3.6);
		\draw [fill=gray1] (4.5,3.2) rectangle (5,3.6);
		\draw [fill=gray1] (5,3.2) rectangle (5.5,3.6);
		
		\draw [fill=gray1] (2,3.6) rectangle (2.5,4);
		\draw [fill=gray1] (5.5,3.6) rectangle (6,4);
		
		\draw [fill=gray1] (2.5,4) rectangle (3,4.4);
		\draw [fill=gray1] (3,4) rectangle (3.5,4.4);
		\draw [fill=gray1] (3.5,4) rectangle (4,4.4);
		\draw [fill=gray1] (4,4) rectangle (4.5,4.4);
		\draw [fill=gray1] (4.5,4) rectangle (5,4.4);
		\draw [fill=gray1] (5,4) rectangle (5.5,4.4);
		\draw [fill=gray1] (5.5,4) rectangle (6,4.4);
		\draw [fill=gray1] (6,4) rectangle (6.5,4.4);
		
		\draw [fill=gray1] (4,4.4) rectangle (4.5,4.8);
		\draw [fill=gray1] (4.5,4.4) rectangle (5,4.8);
		\draw [fill=gray1] (5,4.4) rectangle (5.5,4.8);
		
		\end{tikzpicture}
		\caption{\label{fig:frog} Schematic of a discretized operation field and driving cycle. Resistance times are represented by gray shades. (See also Figures~\ref{fig:nox-nedc} and~\ref{fig:nox-random}.)}
	\end{figure}

	\subsubsection*{Discrete optimization formulation}
	
	In the discrete version of problem~\eqref{eq:continuous-opt}, we
	approximate the objective function in a similar fashion by replacing
	again the integral along the curve $\phi$ by a weighted sum. To this
	end, let as above the weight function
	\[
	\Phi \colon [k]\times[k]\to\R_{\geq 0},\ (f,t)\mapsto \Phi_{ft}
	\]
	be the discrete analogue of the curve $\phi$ describing how long
        the curve $\phi$ stays in $\OP_{ft}$.
	Again, $\Phi$ may coincide with $\dcDisc$ from \eqref{eq:dcDisc} or be chosen to optimize fuel consumption
	on a larger part of the operation field.
	The optimization objective is now to find a solution map $\SOL \in
	\Omega_k$ by picking for each stack $S_{ft}$ a single representative
	$d_{ft}= \SOL(f,t)$, so that $\SOL \in\Omega_k$ solves 
	\begin{equation}\label{eq:discrete-opt}
          \minimize_{M_k\in\Omega_k}\ \sum_{(f,t)\in [k]\times [k]} \Phi_{ft}\cdot [M_k(f,t)]_\fuel \
          \text{subject to the emission constraint~\eqref{eq:em-constraint-disc}.}
	\end{equation}
        Note that $M_k\in\Omega_k$ includes the
	discrete drivability constraint~\eqref{eq:driv-constraint-disc}.
	In Section~\ref{sec:discrete-optimization} we will demonstrate how to formulate this problem
	as an integer linear program (ILP).

	\section{Semi-automatic Calibration}
	\label{sec:calibration}
	\noindent
	Our method consists of two parts, one semi-automatic, and one automatic.
	Algorithm~\ref{algo:calibration} is a first rough sketch of the automatic part, which will be detailed in the following.
	As its input it is given the domain $\domain$, equipped with a grid, and an evaluation oracle for the function $F$.
	It returns a solution map $\SOL\colon [k]\times [k] \to \R^m\times
	\R^n$, $\SOL\in \Omega_k$, of representatives of the $k$-operation field.
	In practice, the evaluation oracle for $F$ is given by an engine mounted on a test-bench.
	
	The while loop of Algorithm~\ref{algo:calibration}, whose description makes up the bulk of this section, terminates if a preliminary solution map
	$\SOLprelim \colon [k]\times [k] \to \R^m\times \R^n$ could be extracted from the measured data. 
	This preliminary solution map $\SOLprelim$ conforms to the emission and drivability constraints.
	However, for each stack $S_{ft}$ in $\kOP$, it contains either a data point $d_{ft}$, which will then represent $S_{ft}$, or a placeholder that adds penalties to the total emission; cf.\ Section~\ref{sec:discrete-optimization}.
	Note that the constraint $M_k\in \Omega_k$ in \eqref{eq:discrete-opt}
	requires, that there are no empty stacks, but our integer linear program (ILP)
	formulation of \eqref{eq:discrete-opt} in \ref{sec:discrete-optimization}, which is used by
	Algorithm~\ref{algo:calibration}, is extended such that it can handle empty stacks by using a
	placeholder instead.
	In the second part, the gaps in $\SOLprelim$ (marked with placeholders) are closed via interpolation-guided measurements and a complete solution map
        $\SOL \in \Omega_k$ is returned; cf.\ Section~\ref{sec:hunt-mode}. 
	The penalty values ensure that replacing a placeholder by a measured value can only improve the solution. 
	As the preliminary solution map $\SOLprelim$ conforms to the given emission constrains, so must the complete solution map
	$\SOL$, which is thus a solution to problem \eqref{eq:discrete-opt}.
	
	In the semi-automatic part of our method, engineering knowledge is used to provide a 
	meaningful base set of measurements for Algorithm~\ref{algo:calibration} and thus 
	guide the calibration process.
	During this so-called basic calibration, the first two steps in Algorithm 
	\ref{algo:calibration} are repeated for a preset amount of time with several actuators
	fixed to values that enforce measurements in regions which are known to be critical
	to any calibration effort, regardless of the specific engine. 
	In particular, low torque regions, which have a significant influence on the overall 
	pollutant emission, are focused on hereby.  
	As this semi-automatic part of our method consists of a subset of the steps in Algorithm 
	\ref{algo:calibration}, we did not dedicate a separate section to its description. 
	Instead, we provide a protocol of the actuator settings and their respective purposes
	during basic calibration in Section~\ref{sec:basic-calibration}.
	
	\begin{algorithm}
		\SetKwInOut{Input}{Input}
		\SetKwInOut{Output}{Output}
		
		\Input{\mbox{admissible domain $\domain$}, \mbox{grid $G$,} \mbox{admissible range $\codomain$,} \mbox{data set $\dmeasurements$,} \mbox{precision parameter $k$,} \mbox{evaluation oracle $F$}}
		\Output{solution map $\SOL$, updated data set $\dmeasurements$}
		
		initialization
		
		\While{no ILP-solution $\SOLprelim$ found}
		{
                  iteration step: adds to $\dmeasurements$ 
			
                  data cleaning: reduces $\dmeasurements$ to $\dmeasurements_{\text{red}}$
                  
                  grid refinement
			
                  check integer linear program from Section~\ref{sec:discrete-optimization} for feasibility with $\dmeasurements_{\text{red}}$ as input and extract solution $\SOLprelim$ if it exists
		}
		
		close the gaps in $\SOLprelim$
		
		\Return $(\SOL, \dmeasurements)$
		
		\caption{Engine calibration procedure}
		\label{algo:calibration}
	\end{algorithm}
	
	%
	%
	
	\subsection{Initialization}\label{sec:init}
	The calibration procedure is initialized with the domain $\domain$ equipped with a grid $G$, the set of noncritical target values $\codomain$, a (possibly empty) set of data points $\dmeasurements$ and a precision parameter $k$.
	The function $F$ is given implicitly; for a given setting of the actuators, the sensors yield the respective function value by means of a physical measurement.
	We assume that the data points in $\dmeasurements$ reflect true values of $F$.
	The domain
	\begin{equation}\label{eq:domain-product}
	\domain \ =\ \prod_{i=1}^{m}I_i
	\end{equation}
	is a product of intervals $I_1,I_2,\dots,I_m$, where $I_i$ defines the range of variation of actuator $i$; cf.\ Section~\ref{subsec:boundaries}.
	It may happen that $I_i$ degenerates to a single point.
	Then the corresponding actuator is called \emph{static}, otherwise it is called \emph{dynamic}. 
	
	There are two cases to distinguish.
	In the first case, the set $\dmeasurements$ of data points is empty.
	Then, for each dynamic actuator $i$, the corresponding interval $I_i$ is subdivided into parts of equal length.
	Otherwise, if $\dmeasurements$ is not empty, then we assume that each interval $I_i$ is equipped with its own subdivision.
	In either case, the product form~\eqref{eq:domain-product} induces a grid structure $G$ on the domain.
	If $\dmeasurements$ is empty then this grid $G$ is uniform.
	In the later stages of the optimization, however, $G$ will become more and more non-uniform.
	
	Each measurand has an interval $J_j$ of noncritical values. 
	For example, the temperature of the test engine has to stay within certain bounds to prevent it from being damaged. 
	Then 
	\[
	\codomain \ =\ \prod_{j=1}^{n}J_j\,.
	\]
	In practice, $\codomain$ is given in part by the physical test engine (e.g., the aforementioned temperature limits) as well as by external factors such as government regulations (e.g., emission limits).
	Just like $F$, which is given by the physical test engine, we will assume the set $\codomain$ to stay fixed throughout the whole calibration process. 
	The precision parameter $k$ is dictated by the engine control unit's engine map format. 
	
	During the various steps of the calibration, data points will be added to the set $\dmeasurements$.
	Thus, it may happen that $\dmeasurements$ becomes prohibitively large to perform subsequent steps of the calibration.
	How to weed out less relevant measurements is the subject of Section~\ref{subsec:cleaning} below.

	\subsection{Iteration Step}\label{sec:iteration-step}
	
	The basic iteration step can be subdivided into two phases: The generation of a measurement plan, followed by the actual measurement, which is combined with a refinement of the grid.
	
	\subsubsection*{Generation of the Measurement Plan}
	
	For our given grid $\grid$ and data set $\dmeasurements$, we construct an abstract graph $\graph=\graph(\grid, E)$ as follows.
	The nodes of $\graph$ are the grid boxes determined by $\grid$.
	Two $m$-dimensional grid boxes are joined by an edge if their intersection is a grid box of dimension $m-1$.
	If there are no measurements yet, i.e., if $\dmeasurements$ is empty, then the grid $\grid$ is uniform, and the graph $\graph$ is the dual graph of a cubical cell complex.
	Due to non-uniform refinement, the structure of $\graph$ will become more complicated.
	
	The graph $\graph$ is equipped with nonnegative node and edge weights. For a grid box $B$, we denote by $\#B$ the number of points $(u^q,y^q)\in\dmeasurements$ such that $u^q\in B$. Then the weight $\boxweight$ of a grid box $B$ is chosen as
	\begin{equation}\label{eq:wB}
	\boxweight(B)\ =\ \frac{\vol(B)}{\#B+1} \,,
	\end{equation}
	which we call the \emph{reciprocal data density} of $B$.
	Adding~1 in the denominator prevents division by zero.
	Further, we define the weight $\edgeweight$ of an edge between two adjacent grid boxes $B$ and $B'$ as the \emph{data density} of $B\cup B'$.
	That is,
	\begin{equation}\label{eq:wBB}
	\edgeweight(B,B')\ =\ \frac{\#(B\cup B')}{\vol(B)+\vol(B')} \,.
	\end{equation}
	There is some room for adjusting these weights; the general idea is that the weights should reflect the data densities.
	
	The classical way to measure the engine behavior is to take a measurement only once a 
	steady-state is reached after an adjustment of the actuators. 
	Accordingly, this technique is called steady-state, or stationary measurement.
	A more recent approach is the quasi-stationary measurement, also called sweep-mapping.
	Here, the actuators are varied slowly and continuously according to a ramp function
	in order to save measurement time. 
	The output then follows with a little delay, while measurements are taken at a regular
	frequency; cf.~\cite[p.~93]{Isermann}.
	There are several techniques to bound and even compensate the contouring error
	of quasi-stationary measurements; cf.~\cite[p.~94ff]{Isermann} and Remark~\ref{rem:hysteresis}.  
	Throughout the while loop of Algorithm~\ref{algo:calibration} quasi-stationary
	measurements will be performed exclusively.
	For these the following concept is crucial.
	
	\begin{definition}(Measurement Ramp)
		For two points $u^q$, $u^r\in\domain$, we call the set
		\[
		\Big\{u^q + i\cdot\frac{u^r-u^q}{\ell-1} \suchthat i=0,1,\dots, \ell-1\Big\}
		\]
		the \emph{measurement ramp} from $u^q$ to $u^r$ with $\ell$ measurements.
	\end{definition}
	
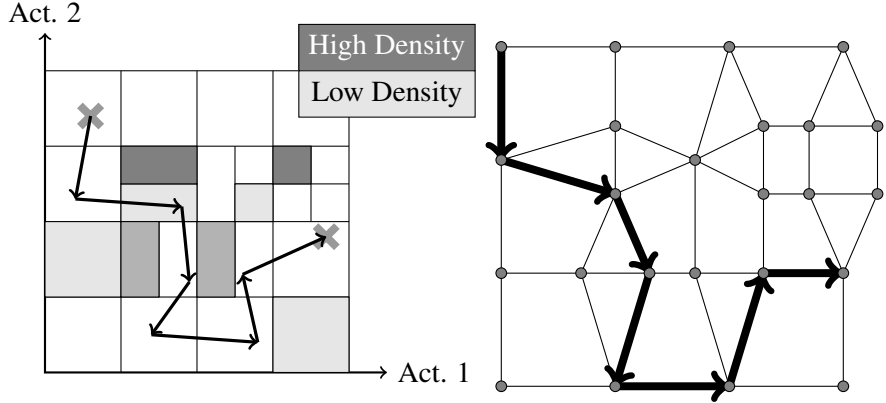
\begin{figure}
		\definecolor{gray1}{gray}{0.9}
		\definecolor{gray2}{gray}{0.7}
		\definecolor{gray3}{gray}{0.5}
		\definecolor{grayCross}{gray}{0.6}
		
		\begin{tikzpicture}[scale=1]
		
		\draw [<->,thick] (0,4.5) node (yaxis) [above] {Act. 2} 
		|- (4.5,0) node (xaxis) [right] {Act. 1};
		
		\draw (1,0) -- (1,4);
		\draw (2,0) -- (2,4);
		\draw (3,0) -- (3,4);
		\draw (4,0) -- (4,4);
		\draw (0,1) -- (4,1);
		\draw (0,2) -- (4,2);
		\draw (0,3) -- (4,3);
		\draw (0,4) -- (4,4);
		\draw (1.5,1) -- (1.5,2);
		\draw (2.5,1) -- (2.5,3);
		\draw (3.5,2) -- (3.5,3);
		\draw (1,2.5) -- (2,2.5);
		\draw (2.5,2.5) -- (4,2.5);
		
		\draw [fill=gray1] (0,1) rectangle (1,2);
		\draw [fill=gray1] (3,0) rectangle (4,1);
		\draw [fill=gray1] (1,2) rectangle (2,2.5);
		\draw [fill=gray1] (2.5,2) rectangle (3,2.5);
		\draw [fill=gray2] (1,1) rectangle (1.5,2);
		\draw [fill=gray2] (2,1) rectangle (2.5,2);
		\draw [fill=gray3] (3,2.5) rectangle (3.5,3);
		\draw [fill=gray3] (1,2.5) rectangle (2,3);
		
		\draw (0.6,3.4) node[cross, line width=3pt,color=grayCross] {};
		\draw (3.7,1.8) node[cross, line width=3pt,color=grayCross] {};
		
		\draw[->, line width=1.2pt] (0.6,3.4) -- (0.4,2.3);
		\draw[->, line width=1.2pt] (0.4,2.3) -- (1.8,2.2);
		\draw[->, line width=1.2pt] (1.8,2.2) -- (1.9,1.2);
		\draw[->, line width=1.2pt] (1.9,1.2) -- (1.4,0.5);
		\draw[->, line width=1.2pt] (1.4,0.5) -- (2.8,0.4);
		\draw[->, line width=1.2pt] (2.8,0.4) -- (2.6,1.3);
		\draw[->, line width=1.2pt] (2.6,1.3) -- (3.7,1.8);
		
		\node[rectangle split, rectangle split parts=2, rectangle split part fill={gray3,gray1}, draw] at (4.5,4) {%
			\textcolor{white}{High Density}
			\nodepart{two}
			Low Density
		};
		\end{tikzpicture}
		\tikzstyle{every node}=[circle, draw, fill=black!50, inner sep=0pt, minimum width=4pt]
		\begin{tikzpicture}[scale=1.5]
		
		\draw (0,0) -- (3,0);
		\draw (0,1) -- (3,1);
		\draw (0,3) -- (3,3);
		\draw (0,0) -- (0,3);
		\draw (0.7,1) -- (1,0);
		\draw (1.3,1) -- (1,0);
		\draw (1.7,1) -- (2,0);
		\draw (2.3,1) -- (2,0);
		\draw (3,1) -- (3,0);
		\draw (0.7,1) -- (1,1.7);
		\draw (1.3,1) -- (1,1.7);
		\draw (1,1.7) -- (1,3);
		\draw (0,2) -- (1,1.7);
		\draw (0,2) -- (1,2.3);
		\draw (1.7,1) -- (1.7,2);
		\draw (1.7,2) -- (2,3);
		\draw (2.3,2.3) -- (2.3,1);
		\draw (2.3,2.3) -- (2,3);
		\draw (1.7,2) -- (1,1.7);
		\draw (1.7,2) -- (1,2.3);
		\draw (1.7,2) -- (2.3,2.3);
		\draw (1.7,2) -- (2.3,1.7);
		\draw (3.3,2.3) -- (2.3,2.3);
		\draw (3.3,1.7) -- (2.3,1.7);
		\draw (3.3,2.3) -- (3.3,1.7);
		\draw (2.7,2.3) -- (2.7,1.7);
		\draw (2.7,1.7) -- (3,1);
		\draw (3.3,1.7) -- (3,1);
		\draw (2.7,2.3) -- (3,3);
		\draw (3.3,2.3) -- (3,3);
		
		\draw[->,line width=3pt] (0,3) -- (0,2);
		\draw[->,line width=3pt] (0,2) -- (1,1.7);
		\draw[->,line width=3pt] (1,1.7) -- (1.3,1);
		\draw[->,line width=3pt] (1.3,1) -- (1,0);
		\draw[->,line width=3pt] (1,0) -- (2,0);
		\draw[->,line width=3pt] (2,0) -- (2.3,1);
		\draw[->,line width=3pt] (2.3,1) -- (3,1);
		
		\node at (0,0) {};
		\node at (1,0) {};
		\node at (2,0) {};
		\node at (3,0) {};
		\node at (0,1) {};
		\node at (0.7,1) {};
		\node at (1.3,1) {};
		\node at (1.7,1) {};
		\node at (2.3,1) {};
		\node at (3,1) {};
		\node at (0,2) {};
		\node at (1,1.7) {};
		\node at (1,2.3) {};
		\node at (1.7,2) {};
		\node at (2.3,1.7) {};
		\node at (2.3,2.3) {};
		\node at (2.7,1.7) {};
		\node at (2.7,2.3) {};
		\node at (3.3,1.7) {};
		\node at (3.3,2.3) {};
		\node at (0,3) {};
		\node at (1,3) {};
		\node at (2,3) {};
		\node at (3,3) {};
		\end{tikzpicture}
		\caption{\label{fig:path}The routing during measurement planning prefers regions of low data density. On the right: Corresponding path in the induced abstract graph.}
	\end{figure}

	Zero-entries of $u^r-u^q$ correspond to actuators which are locally static.
	Let $\numberF\define\abs{\dmeasurements}$ and suppose that there is an admissible point $u=u^\numberF\in\domain$ where the last measurement took place.
	If this does not exist, we choose $u$ uniformly at random in the domain $\domain$.
	Let $B$ be the grid box containing~$u$.
	We may assume that $B$ is unique, since $u$ has been constructed in a randomized fashion.
	The two steps of the generation of the measurement plan are as follows:
	\begin{enumerate}[I.]
		\item\label{it:random-grid-box} \emph{Random grid box:}
		Pick a grid box $B'$ at random with probability 
		\[
		\frac{\boxweight(B')}{W} \,,
		\] 
		where $W=\sum_{B\in G} \boxweight(B)$ is the sum of the reciprocal data densities of all boxes.
		\item \emph{Measurement path:}
		Determine a shortest path $B=B_0,B_1,\dots,B_s=B'$ in $\graph$ from $B$ to $B'$ using Dijkstra's algorithm with respect to the weights
                $\edgeweight$ in \eqref{eq:wBB}; cf.\ Figure~\ref{fig:path} and~\cite[\S2.2]{Cook1998}.
		In each box $B_q$ for $q \in [s]$, pick a point $u^{\numberF+q}$ uniformly at random.
		For $1\leq j\leq s$, connect the points $u^\numberF,u^{\numberF+1},\dots,u^{\numberF+s}$ by $s$ measurement ramps
		with~$\ell_j$ measurements each.
	\end{enumerate}
	This will result in a total of up to $\sum_{j=1}^s\ell_j$ new measurements to be added to the set $\dmeasurements$. 
	In our setting, the actuators are varied at a constant speed on each measurement ramp. 
	This speed is adjusted so that at least one actuator is varied at its maximal variation speed;
	the values for the maximal actuator variation speeds are listed in Table~\ref{tab:actuators}. 
	The measurement frequency throughout the whole calibration process is set to one measurement
	per second. 
	As a consequence, both the length and the orientation of the measurement ramps determine the numbers $\ell_j$.
	
	\subsubsection*{Measurement}
	
	The actual measurement is combined with an adaptive refinement of the grid.
	The goal of the iteration is to fill the solution map in a way that the remaining small holes can be closed by extrapolation of the surrounding data. 
	The next step is then:	
	\begin{enumerate}[I.]
		\setcounter{enumi}{2}
		\item \emph{Measurement:} 
		The actuator settings are varied continuously along the path prescribed by the measurement ramps.
		This results in a linear ordering of the measurements.
	\end{enumerate}
	
	For the recording of the $\ell_j$ measurements per ramp, several additional aspects have to be taken into account, as explained below.
	It is important not to store every measurement, since we do not want to store redundant data.
	Thus, we restrict our attention to \emph{relevant value variations}, i.e., we only store the measurement $F(u^q)$ at a point $u^q$ if the measurement is sufficiently different.
	Consequently, fewer than $\sum_{j=1}^s \ell_j$ data points may be added to the set $\dmeasurements$.
	Here one can choose among various meaningful distance functions.
	However, we also want to establish a \emph{minimal measurement frequency}.
	That is, if we omitted too many subsequent measurements due to the previous rule, then we store the measurement nonetheless.
	
	Moreover, it may happen that a measured value lies outside the admissible range $\codomain$, i.e., it violates one or more restrictions.
	In order to \emph{exclude critical values}, the entire measurement is disrupted, and we continue from scratch at~\eqref{it:random-grid-box} with the last valid measurement.
	Special care is needed for observables with a pronounced latency.
	The development of these values is extrapolated, and the measurement is rerouted already if the measured values get sufficiently close to the boundary of $\codomain$.

	\subsection{Data Cleaning}\label{sec:data-cleaning}
	\label{subsec:cleaning}
	
	Assume that the measurement phase of the iteration step is complete, i.e., a new set
	\[
	\dmeasurements\ =\ \{(u^1,y^1),(u^2,y^2),\dots,(u^{\numberF},y^{\numberF})\}
	\]
	is available. In Section~\ref{sec:grid-refinement}, we will compute local polynomial fits to the function $F$ using data points in $\dmeasurements$ as interpolation points. At points in which the interpolation is not good enough, the grid is refined further. 
	A fundamental necessity for the good fit of a polynomial approximation is that the interpolation points from which it is generated resemble a random point cloud. This is not the case for $\dmeasurements$, since all points lie on a piecewise linear path given by the sequence of measurement ramps. We thus need to extract a subset $\dmeasurements_{\text{red}}$ from $\dmeasurements$ that is sufficiently generic.
	
	Furthermore, we want to translate the discrete version~\eqref{eq:discrete-opt} of our optimization problem into an integer linear program which picks the right kind of measurements from which we can then obtain our engine maps; cf.\ Section~\ref{sec:discrete-optimization} below.
	However, our measurements need to be preprocessed in order to make such an approach feasible by significantly reducing the amount of data without affecting the accuracy.
	
	As it turns out, we can achieve both goals with a single method, which we call the \emph{adaptive space compressor}.
	The idea is the following.
	For the set $\widehat U = \{ u^1,u^2,\dots,u^\numberF \}$, which is the projection of the set $\dmeasurements$ onto the admissible domain $\domain$, we define a corresponding \emph{threshold graph}.
	Its nodes are the points $u^q$, and there is an edge between $u^q$ and $u^r$ if the Euclidean distance $\norm{u^q-u^r}_2$ is below a certain threshold.
	This threshold depends on the sizes of the two grid boxes which contain $u^q$ and $u^r$.
	Then we iteratively remove measurements which have the maximal node degree in the threshold graph until only isolated nodes remain.
	This \emph{reduced} set is denoted as $\widehat U_{\text{red}}$, which yields
	\[
	\dmeasurements_{\text{red}}\ =\ \bigl\{ (u,y) \in \dmeasurements \suchthat u \in \widehat U_{\text{red}} \bigr\}\,.
	\]
	The set $\dmeasurements_{\text{red}}$ has significantly fewer data points compared to $\dmeasurements$. 
	As the threshold for the connection of $u^q$ and $u^r$ by an edge depends on the size of the
	grid boxes containing them, the distribution of the remaining data points reflects the structure
	of the grid $G$.
	Due to the criteria for grid refinement, which will be introduced in (V) and (VI) below, the structure of the grid $G$ reflects the behavior of $F$.
	In particular, this means that the data density is higher in regions where the behavior of $F$ is hard to reconstruct by interpolations. 
	Thus, the adaptive space compressor filters out the relevant information gained by the preceding measurements. 
	Moreover, it breaks up the piecewise linear structure of the data point distribution, 
	leading to local (pseudo-)randomness which is required by the polynomial interpolation. 
	
	Note that this reduction does not affect the set $\dmeasurements$, but we rather explicitly keep $\dmeasurements_{\text{red}}$ as a subset.
	In this way, all our measurements are taken into account for computing the weights according to~\eqref{eq:wB} and~\eqref{eq:wBB} in (I) in the next round of measurements.	
	
	\begin{remark}\label{rem:hysteresis}
		Since the points $u^i$ are picked in a randomized fashion, the piecewise differentiable map $F$ is differentiable at $u^i$ almost surely. 
		Thus, we may assume that every $u^i$ lies at the center of an open ball on which $F$ is smooth. 
		The grid structure is refined and thus the measurement density increased in areas which are not sufficiently smooth; cf.\ Section~\ref{sec:grid-refinement}. 
		Consequently, in an individual grid box one can, with high probability, cluster the measured points according to smooth patches of $F$. 
		As a result, for such clusters the local time constants are identifiable and the techniques for reversing the contouring error by dynamic correction in~\cite[p.~95f]{Isermann} can be applied successfully, leading to stored measured values
		that closely approximate the steady state values of the engine. 		
	\end{remark}
	
	\subsection{Grid Refinement}\label{sec:grid-refinement}
	
	The purpose of refining the grid $\grid$ is to accumulate sufficiently many \enquote{meaningful} data points in the set $\dmeasurements_{\text{red}}$, such that the $k^2$ representatives for the solution map can be extracted while observing the emission and drivability constraints.
	Recall that in our application we have $m\geq 8$ actuators and $n\geq 14$ measurands. This makes computing on a uniform grid infeasible.
	The price to pay is that it is slightly more involved to determine points $u^q$ in the admissible domain $\domain$ such that the data point $(u^q,y^q)$ adds significant information to our data set $\dmeasurements_{\text{red}}$ and is thus stored.
	We proceed with the following step:
	\begin{enumerate}[I.]
		\setcounter{enumi}{3}
		\item \emph{Computation of local fit:}
		We employ a Newton interpolation approach, which requires only a partial recomputation of the interpolation polynomial if some interpolation points are exchanged. Each component function $[F]_i: \Rm\rightarrow \R$, where $i\in[n]$, has to be interpolated by a Newton polynomial $L_d[F]_i$, where $d$ is the approximation order. 
		Assuming the $(d+1)$-times continuous differentiability of $[F]_i$ at $u$, a $d$-th order polynomial fit requires $N \define {{m+d}\choose{m}}$ interpolation points $\{u^{q_1}, u^{q_2}, \dots, u^{q_N}\}\subseteq \umeasurements$.
		These $N$ points have to satisfy certain spatial relations for the approximation properties of the interpolation polynomials to hold; the interpolation problem is then called \emph{poised}.
		The formulation of the latter conditions is somewhat technical.
		We thus omit the details and refer to the literature instead; cf., e.g.,~\cite{Schueler2000}.
		Our data cleaning method in Section~\ref{sec:data-cleaning} ensures poisedness (with high probability).
		
		If an interpolation polynomial with $N$ interpolation points which have the required geometrical structure fails to give a $d$-th order approximation of $[F]_i$ at $u$, this indicates insufficient smoothness of $[F]_i$ at $u$, which we will use as a criterion for a local grid 
		refinement.
		We aim at a second order fit, i.e., $d=2$, which requires
		\[
		N\ =\ {{m+2}\choose{m}}\ =\ \frac{m\,(m+1)}{2} 
		\]
		interpolation points. Our intention is to compare the measured value $y^q\in\codomain$ at a point $u^q\in\domain$, where $(u^q,y^q)\in\dmeasurements_{\text{red}}$, with a polynomial interpolation of $F$ in $u^q$ using elements of $\dmeasurements_{\text{red}}$ as interpolation points.
		The error of such an approximation is minimized if we use the $N$ closest neighbors $\{u^{q_1}, u^{q_2}, \dots, u^{q_N}\}$ of $u^q$ in $\domain$ as interpolation points, such that 
		\[
		\{(u^{q_1},y^{q_1}), (u^{q_1},y^{q_2}), \dots, (u^{q_N},y^{q_N})\}\ \subseteq\ \dmeasurements_{\text{red}}\,.
		\] 
	\end{enumerate}
	
	Concerning the implementation of the polynomial interpolation: As can be seen in the experimental section below, we will usually have less than $100,000$ points in $\dmeasurements_{\text{red}}$. This makes the following brute-force approach feasible.
	Let $\numberF_{\text{red}} \define \card{\dmeasurements_{\text{red}}}$ and assume for simplicity that $q=1$.
	Now calculate the squared Euclidean distances of $u^2, u^3, \dots, u^{\numberF_{\text{red}}}$ to $u^1$. This has a cost of roughly $3\cdot m\cdot \numberF_{\text{red}}$ elementary arithmetic operations. Store these values in the first row of a $(2\times [\numberF_{\text{red}}-1])$-array and the indices of the corresponding points $u^r$ in the second row ($\approx 2\cdot \numberF_{\text{red}}$ writes). Then determine the $N$ smallest entries of the first row and return the corresponding second row entries. 
	
	A simple in-place algorithm to accomplish this is the following: First, determine the $N$-th smallest element of the first row, which costs one sweep of the array; cf.~\cite[p.~183 ff]{Cormen:2009:IAT:1614191}. Second, sort all columns to the front of the array whose first row entry is smaller than or equal to the $N$-th smallest first row entry. This costs another sweep of the array. Finally, return the first $N$ second row entries.
	
	Repeating this procedure for all $\numberF_{\text{red}}$ points in $\dmeasurements_{\text{red}}$, we arrive at an approximate cost of $(3\cdot m+4)\cdot \numberF_{\text{red}}^2 \in \mathcal O(m\cdot \numberF_{\text{red}}^2)$ elementary operations.       
	In our setting, for up to $100,000$ points and $8$ dynamic actuators, this results in a total of about $280$ billion elementary operations, which is a fairly insignificant task for modern computers.    
	For a general discussion on the nearest neighbor search in high dimensions we refer to the survey~\cite{AndoniIndyk:2018}.
	The next step is the following:  
	\begin{enumerate}[I.]
		\setcounter{enumi}{4}
		\item \emph{Symmetric grid refinement:}
		To decide whether the grid needs to be refined, we compare the interpolated value, say $\tilde F(u^q)$, with $y^q=F(u^q)$.
		If the deviation exceeds a threshold, then the box $B$ which contains~$u^q$ is split into $2^m$ smaller grid boxes, symmetrically in all coordinate directions.
		
		The quality of the local fit depends on the differentiability properties of the approximated function. Hence, the symmetric grid refinement increases the measurement density in areas of potential nondifferentiability. 
	\end{enumerate}
	In practice, the static actuators can be ignored for the computation of the local fit and the symmetric grid refinement.
	
	\subsubsection*{Making the Grid Nonuniform}
	
	Recall that ultimately we want to solve the  discrete version~\eqref{eq:discrete-opt} of the constrained optimization problem~\eqref{eq:continuous-opt}, i.e., we want to minimize the fuel consumption subject to the drivability and emission constraints.
	For the iteration step, the drivability constraints are irrelevant.
	Therefore we focus on the subset $\emission$ of measurands corresponding to the emissions.
	The final output of the calibration will be a solution map
	$\SOL\in\Omega_k$ represented by $k^2$ data points $d_{ft}=\SOL(f,t)$, where $f,t\in[k]$, which cover the $k$-operation field,  while satisfying the emission constraints~\eqref{eq:em-constraint-disc}. 
	In practice, it is a major challenge to find sufficiently many points which satisfy the conditions~\eqref{eq:em-constraint-disc} imposed by emission control.
	This leads us to a second type of grid refinement.
	
	\begin{enumerate}[I.]
		\setcounter{enumi}{5}
		\item \emph{Asymmetric grid refinement:}
		Consider the point $u^\numberF\in\domain$ at which the last measurement was performed.
		Let $f$ and $t$ be indices such that $(u^\numberF_{\frequency},y^\numberF_{\torque})\in\OP_{ft}$, where $y^\numberF=F(u^\numberF)$.
		If, for any $p\in\emission$, we have
		\[
		[F(u^\numberF)]_p\ < \ \min_{(u^q,y^q)\in S_{ft}}([F(u^q)]_p)\,,
		\]
		i.e., if the emission measurement $[F(u^\numberF)]_p$ is lower than any other value on the corresponding stack $S_{ft}\subset\dmeasurements$, a \emph{cross-measurement} is performed.
		To this end, each actuator is varied individually to determine the direction with the biggest impact on $[F]_p$.  
		Afterwards, the grid box $B$ containing $u^\numberF$ is split into two congruent (sub-)boxes
		along the axis corresponding to the actuator whose variation has the biggest impact on $[F]_p$.
	\end{enumerate}
	Both types of grid refinement do not add data points to $\dmeasurements_{\text{red}}$ directly.
	However, the grid refinements increase the probability for picking one of the subboxes in (I).
	This way one can hope to find data points near $(u^q,y^q)$ that add relevant information about $F$ and near $(u^\numberF,y^\numberF)$ with better emission values than those currently stored in $\dmeasurements$.

	\subsection{An Integer Linear Program}\label{sec:discrete-optimization}
	
	As advertised in Section~\ref{sec:discretization}, we will now detail how to 
	formulate the discrete optimization problem~\eqref{eq:discrete-opt} as an 
	integer linear program (ILP) in order derive a preliminary solution map
	$\SOLprelim\colon [k]\times [k] \to \R^m \times \R^n$ from the reduced set of measurements $\dmeasurements_{\text{red}}$ that satisfies the objectives listed in Section~\ref{sec:discretization}. 
	This ILP-solution $\SOLprelim$ is not necessarily a complete solution map $\SOL$, as it may 
	contain some placeholders, which must be replaced in a subsequent interpolation step; 
	cf.\ Section~\ref{sec:hunt-mode}. 
	For an introduction to the solution of ILPs, see~\cite{Cook1998},~\cite{Schrijver1998}, and~\cite{Sierksma2001}.
	
	For a data point $d^q$, let $S^q$ be the stack that contains it.
	As before, let \enquote{$\fuel$} be the data point index corresponding to the injected fuel quantity.
	Then we call the weight 
	\[
	\Phi^q \define
	-\frac{\min\bigl\{[d^r]_\fuel \suchthat d^r\in S^q \bigr\} }{ [d^q]_\fuel}
	\]
	the \emph{prey value} of $d^q$. 
	The prey value of $d^q$ is the negative of the quotient of the minimal fuel 
	consumption among all data points in $S^q$, the stack containing $d^q$, over the fuel consumption at $d^q$. 
	As such, the absolute value of any prey value is smaller than or equal to $1$. 
	Ideally, a prey value equals $-1$.  
	
	Further, for each index $q\in[\numberF_{\text{red}}]$, where $\numberF_{\text{red}} \define \card{\dmeasurements_{\text{red}}}$, we introduce a binary decision variable $s^q\in\{0,1\}$, which indicates whether the data point $d^q$ is part of the ILP-solution $\SOLprelim$.
	The objective function of our integer linear program can now be written as
	
	\begin{equation}\label{eq:opt}
	\minimize\  \sum_{q=1}^{\numberF_{\text{red}}} \Phi^q \, s^q\,,
	\end{equation} 
	while conforming to the constraints~\eqref{eq:stack-constraint},~\eqref{eq:drivability-constraints}, and~\eqref{eq:integral-emission-constraint}, 
	which are described in detail below. 
	It was already noted in Section~\ref{sec:continuous-opt} that there exists a wide array 
	of meaningful weight functions $\Phi$ (and their continuous counter-parts $\phi$).
	The above choice for $\Phi$, the prey values, ensure that the ILP-solver picks for 
	all operation points the data point with the least possible fuel consumption among 
	all feasible choices. 
	
	\subsubsection*{The stack constraint}
	The solution $\SOLprelim$ of our ILP should contain at most one element from each stack $S_{ft}$.
	This condition is reflected in the \emph{stack constraint}
	\begin{align}\label{eq:stack-constraint}
	s_{ft}\ +\ \sum_{d^q\in S_{ft}}s^q\ =\ 1 \quad \text{for all }\  (f,t)\in[k]\times[k]\,,
	\end{align}
	where $s_{ft}$ is a \emph{stack decision variable}; it is $1$ if it is not possible to choose a data point for stack $S_{ft}$ and $0$ otherwise.
	The nonzero $s_{ft}$ are the placeholders we mentioned above. 
	They will be used in~\eqref{eq:integral-emission-constraint} to assign the  
	penalty values to stacks that contribute no data point to the ILP-solution $\SOLprelim$. 
	
	\subsubsection*{Formalizing the drivability constraint} 
	To avoid engine damage, the discrete drivability constraint \eqref{eq:driv-constraint-disc}
        ensures, that the variation speeds of all $m$ actuators are constrained individually by nonnegative constants $\Delta_\actuator$ for $\actuator\in[m]$.
	To write this condition as a constraint in an ILP, let $d^q\in S_{ft}$, and $d^r$ be a data point in either neighboring stack $S_{f\pm 1,t}$ or $S_{f,t\pm 1}$.
	Then the \emph{drivability constraints} \eqref{eq:driv-constraint-disc}
	can be expressed as
	\begin{align}\label{eq:drivability-constraints}
	s^q\ +\ s^r\ \leq\  1\qquad\text{if}\quad \bigl\lvert [d^q]_\actuator-[d^r]_\actuator \bigr\rvert\ \geq\ \Delta_\actuator\quad \text{ for any }\ \actuator\in[m]\,.
	\end{align}
	Note that this is a secant constraint, since the data points compared are contained in neighboring stacks. 
	
	\subsubsection*{Formalizing the emission constraint} 
	The emission condition~\eqref{eq:em-constraint-disc} describes the upper bound of several emission test cycles, e.g., maximal \NX production.
	Additionally, the current output of a pollutant $p\in\emission$ is restricted by the boundaries of the respective interval of noncritical values
	\[
	J_p\ =\ [\underline{e}_p,\, \overline{e}_p]\,.
	\]
	In the following, the weights $\omega_{ft}$ represent the mean resistance time of the rectangles $\OP_{ft}$ in the given test cycle.
	For simplicity of notation, the mean resistance time (cf.\ Section~\ref{sec:discretization}) on the rectangle corresponding to the stack $S^q$ containing $d^q$ is denoted by $\omega^q$.
	Then the \emph{emission constraint} can be formulated as follows: 
	\begin{align}\label{eq:integral-emission-constraint}
	\sum_{q=1}^{\numberF_{\text{red}}} \omega^q\, [d^q]_p\, s^q\ +\ \sum_{(f,t)\in [k]\times [k]} \omega_{ft}\, \overline{e}_p\, s_{ft}\ \leq\ \emconstraint_p\qquad\text{for all }\ p\in\emission\,.
	\end{align}
	The second summand of the left hand side of the inequality serves as the above-mentioned 
	penalty term which compensates for stacks that contribute no data point to the ILP-solution $\SOLprelim$. 
	
	\begin{remark}
          We want to determine the size of our ILP.
          There are $k^2$ equalities arising from the stack constraints.
          Further, it has $m\cdot 2\cdot k\cdot(k-1)$ drivability constraints and
          $\vert E\vert$ linear emission inequalities.
		There are $k^2$ stack decision variables and $\numberF_{\text{red}}$ decision variables
		$s^q$, one for each data point in $\dmeasurements_{\text{red}}$.
		In the situation of the simulation data presented below, where $k=16$, and $m=8$, this results
		in $256$ equalities and $3840+8$ inequalities which constrain the ILP-solution $\SOLprelim$.
		Moreover, there are $256$ stack decision variables. 
		The size of $\dmeasurements_{\text{red}}$ naturally varies throughout the calibration process, 
		and may range between $10.000$ and $100.000$.
	\end{remark}

	\subsection{Closing the Gaps in $\SOLprelim$}\label{sec:hunt-mode}
	
	The calibration algorithm traces the behavior of the map $F$, given by the physical test engine, by constructing a sequence of (one-dimensional) measurement ramps through the domain $\domain$ which is high-dimensional, as is the range of $F$. 
	Naturally, this approach cannot produce a sufficient coverage of $\domain$.  
	In particular, often several stacks $S_{ft}$ of the $k$-operation field will contain no 
	data point that contributes to the solution $\SOLprelim$ of the integer linear program. 
	
	For each such non-contributing stack $S_{ft}$ we proceed as follows. 
	First, we construct a local model of $F$ as follows.
	If $S_{ft}$ lies in the interior of the operation field, we
	pick $N=m(m+1)/2$ data points
	\[
	\{(u^{q_1},y^{q_1}), (u^{q_2}, y^{q_2}),\dots,(u^{q_N},y^{q_N})\}
	\] 
	from neighboring stacks and compute a Newton type polynomial that interpolates
	the points $\{u^{q_1}, u^{q_2}, \dots,u^{q_N}\}$ as done in the local fit step IV.
	If an empty stack $S_{ft}$ lies on the boundary of the operation field, then the 
	local model is computed by calculating secants of neighboring data points 
	in the nearest stacks and extending them linearly. 
	
	Second, we find, e.g., via Newton's method, some point $\tilde u\in\domain$ such that under our local model of $F$
	we have $(\tilde u_\frequency,\tilde y_\torque)\in S_{ft}$.  
	In a ball about $\tilde u$ we perform randomized measurements, e.g., using a normal distribution
	centered at $\tilde u$, until we find a point $u$ with
        $(u_\frequency, y_\torque)\in S_{ft}$ and which satisfies the drivability constraints.
	Due to the continuity of the manifold map there must exist a neighborhood of such points. 
	The data point $d_{ft}\define (u,y)$ is then added to the solution as the representative of the stack $S_{ft}$. 
	Our way of picking the $d_{ft}$ ensures that the completed solution map conforms to the drivability constraint and we thus have $\SOL\in\Omega_k$.
	
	For pollutant $p$ the penalty value is $\overline{e}_p$, the upper bound of the interval 
	$J_p = [\underline{e}_p,\, \overline{e}_p]$ of noncritical values. 
	Hence, any recorded value of pollutant $p$, by which the penalty value is replaced, is 
	smaller than or equal to $\overline{e}_p$. 
	As a consequence, the total emission of a complete solution map $\SOL$ cannot exceed that of
	$\SOLprelim$ and must thus conform to the emission constraint if the preliminary solution does. 
	Since we also have $\SOL\in\Omega_k$, as noted above, the so-completed solution map does indeed solve problem \eqref{eq:discrete-opt}.
	
	\section{AVL Engine Model}\label{sec:avl}
	
	For our experiments in Section~\ref{sec:experiments}, we replace the test-bench with a model of a diesel engine with turbo charger, pilot injection and variable turbine geometry.
	The latter has been developed in cooperation with AVL GmbH and is based on measurements on a compression ignition/diesel engine.
	Below we will give a brief overview of the effects of the $8$ actuators and $18$ measurands that are simulated, thus indicating the scope of the simulation. 
	For a detailed description of the AVL model's derivation, see~\cite[p.~69ff]{Burggraf} and~\cite{AVL}.
	
	\subsection{Actuators of the AVL Engine Model}\label{sec:actuators}
	
	\paragraph{\it Revolution frequency of the crankshaft (\RF)}
	The crankshaft converts the reciprocating motion of the cylinders into a rotational motion.
	In modern four-stroke engines every cylinder fires once for every two revolutions of the crankshaft.
	The revolution frequency is a controlled actuator, as discussed in Section~\ref{sec:math}.
	It stands out among the other actuators since it provides one coordinate axis of the operation field.
	
	\paragraph{\it Injected fuel quantity (\IF)} In contrast to a spark-ignited engine, the injected amount of fuel is the most important actuator. More precisely, the injection process is crucial in the application process of diesel engines. The injection process is given by several pre/pilot-injections, a main injection and post-injections. Typically, the engine torque is mainly determined by the main-injection. This actuator defines the total amount of fuel per cycle. In this simple model, the injected fuel volume is divided into one pilot and the main injection.
	
	\paragraph{\it Pressure in the common rail system (\RP)} In contrast to solenoid-controlled unit injector elements, the pressure is generated by a central fuel and high pressure pump. The fuel injectors are opened and closed by piezo elements. 
	
	\paragraph{\it Air filling (\AF)} Similar to a spark ignited engine, an air valve controls the amount of air which contributes to the combustion process. In this model the amount of air is given directly in mass per piston stroke. 
	
	\paragraph{\it Turbine geometry (\TG)} Modern turbochargers do not have a static turbine geometry. Variable-turbine-geometry turbochargers are able to tune the angle of the turbine blades in order to increase the amount of boost. Alternative setups are given by static turbines with waste gates. Waste gates are applied to reduce the amount of exhaust gas that accelerates the turbine, so the amount of boost can be controlled. In our model the geometry is given as a value between $30$ and $85$. 
	
	\paragraph{\it Main timing (\MT)} The main timing is comparable to the spark timing of Otto-engines. It defines the start timing of chemical reactions in the crank\-shaft angle of the main injection. Similar to the Otto engine, the pressure rise is delayed by the ignition delay. 
	
	\paragraph{\it Pilot injection (\PI) and pilot timing (\PT)} Pilot injection works in tandem with pilot timing to achieve a complete burning of the fuel, which in turn also drastically reduces the emission of \NX gases. The pilot injection increases the temperature of the combustion chamber, thus when the main injection occurs the fuel is sent into a chamber which already is at a higher temperature than its autoignition point. This especially facilitates the fuel burning at lower speeds.
	
	\renewcommand{\arraystretch}{1.3}
	\begin{table}	
		\caption{\label{tab:actuators} Actuator-intervals during various calibration runs. 
			$\Delta^{\max}$ denotes the maximal variation speed of the actuator in respective units per second. LT/HT = low/high torque, FR = free variation of torque; cf.\ Sections~\ref{sec:basic-calibration} and~\ref{sec:full-calibration}. }
		\small
		\centering
		\begin{tabular*}{\textwidth}{@{}l@{\extracolsep{\fill}}llllll@{}}\toprule
			Act. & Unit & $\Delta^{\max}$ & Basic, LT & Basic, HT & Basic, FR & Full Calibration \\ \midrule
			\RF & $\frac {1}{\text{min}}$ &10& 1000--2600 & 1000--2600 & 1000--2600 & 1000--2600 \\
			\IF & $\frac{\text{mm}^3}{\text{cycle}}$ &0.1& 6--10 & 50--60 & 6--60 & 6--60 \\
			\RP & $\operatorname{hPa}$ &100& 295677 / 405677 & 295677 & 295677 & 295677--1126537 \\
			\AF & $\frac{\text{mg}}{\text{stroke}}$ &5& 300 & 300 & 300 & 275--991 \\
			\TG & $\operatorname{int}$ &1& 30 & 30 & 30 & 30--85 \\
			\MT & $ ^\circ \operatorname{CA}$ &0.2& 0 / 10 & 0 & 0 & 0--10 \\
			\PI & $\frac{\text{mm}^3}{\text{cycle}}$ &0& 1 & 1 & 1 & 1 \\
			\PT & $\mu s$ &10& 1540 & 1540 & 1540 & 1540--2565	\\\bottomrule
		\end{tabular*}
	\end{table}
	
	\subsection{Measurands of the AVL Engine Model}
	\paragraph{\it Torque} The produced torque of the engine. The revolution frequency and the engine torque define the power level of the ICE. 
	It stands out among the other measurands since it provides the second coordinate axis of the operation field.
	
	\paragraph{\it Fuel mass flow} The amount of fuel which is delivered to the cylinder per hour.
	
	\paragraph{\it Carbon monoxide, hydrocarbons, nitrogen oxides, soot}  The momentary and integral exhaust emission of pollutants, given in parts per million and emission per hour. 
	
	\paragraph{\it Indicated mean pressure} The average pressure over a cycle in the cylinder.
	
	\paragraph{\it Lambda} Lambda displays the chemical partitioning of the fuel. Spark-ignited engines run on lambda around $1$. The combustion limits are 0.6 and 1.6. The lambda value for compression ignition engines is much higher, up to 20.
	
	\paragraph{\it Manifold pressure} The manifold pressure measures the absolute pressure in front of the intake channel. The pressure depends on the absolute environmental pressure and the boost level of the turbocharging system. 
	
	\paragraph{\it Boost pressure} The boost pressure is created by the turbocharging unit. It depends on the turbine geometry setting and the current combustion behavior.
	
	\paragraph{\it Maximal cylinder pressure} The maximal point of the cylinder pressure sequence. Every engine has a specified maximal cylinder pressure, in order to avoid damaging of the devices. Typically, the maximal pressure is about $160$ bar.
	
	\paragraph{\it Manifold temperature} The manifold temperature measures the temperature in the intake channel.
	
	\paragraph{\it Critical temperature} The critical temperature is defined as the temperature of the burn zone when the exhaust valves open. In case of bad timings, the fuel has not been consumed completely, which results in the release of flames to the exhaust manifold, catalysts and turbocharger. It indicates damages to sensitive parts of the engine setup. 
	
	\paragraph{\it Specific fuel consumption} The current power level of the engine over the current fuel consumption.

	\section{Optimizing the AVL Engine Model}\label{sec:experiments}
	
	As explained in Section \ref{sec:calibration}, we subdivide the engine optimization into two phases.
	During the first phase, which we call \emph{basic calibration}, measurements are taken for a prescribed
	amount of time in regions which are known to be critical to any calibration effort, regardless
	of the specific engine. 
	The basic calibration is concluded by passing the obtained data to Algorithm~\ref{algo:calibration}, which is run without the emission constraint~\eqref{eq:integral-emission-constraint} in the ILP. 
	If constraint~\eqref{eq:integral-emission-constraint} is omitted, there always exists a feasible solution
	$\SOLprelim$ to the ILP, which may contain several gaps, though, due to the drivability constraint and empty stacks. 
	Hence, the while loop exits immediately and the gaps in $\SOLprelim$ are closed. 
	
	In the second phase, which we call \emph{full calibration}, a solution map conforming to all constraints is obtained by building on and refining the preliminary solution map of the first phase. 
	During full calibration, Algorithm~\ref{algo:calibration} runs automatically and without time limit until it returns a complete solution map $\SOL$.
	
	One could, in theory, omit the first phase and turn the calibration procedure into a fully
	automatic one by letting all actuators range freely on their whole domain right from the start. 
	However, carefully applying engineering-knowledge in the first phase by bounding some 
	actuators, fixing others, and ignoring the very restrictive emission constraints results
	in the quick gathering of a meaningful base set of data which then merely has to be refined in
	the sequel.
	In our experience, this approach considerably reduces the number of measurements required.

        Throughout this section the precision parameter $k$ is set to 16, and this yields 256 representative data points for the solution map $\SOL$.
        This determines the shapes of Figures~\ref{fig:af-comp}ff.
        
	\begin{remark}
		In the test bench scenario, measurements are taken, though not necessarily stored, in regular time intervals. 
		A realistic frequency is one measurement per second. 
		We take this frequency as the basis of our translation of the number of measurements in the
		simulation into real-world time.
		Below, we state the costs of our method in real-world time, as all computation times occurring
		during the different steps of Algorithm~\ref{algo:calibration} are negligible in comparison to
		the days, or even weeks, that the physical experiments on a test bench can take.
	\end{remark}
	
	\subsection{Basic Calibration}\label{sec:basic-calibration} During basic calibration, only two actuators are dynamic, \IF and \RF. The other six actuators are static throughout. 
	The phase is subdivided into three runs of several hours. Here a \emph{run} means the following: 
	The intervals of the dynamic actuator \IF and the values of the static actuators are reset. Then the algorithm steps in Sections~\ref{sec:iteration-step}--\ref{sec:grid-refinement} are repeated for a preset amount of time.
	The revolution frequency is set to vary over its full range of 1000--2600 revolutions per minute throughout all three runs.
	For a full account of the applied settings, see Table~\ref{tab:actuators}.
	
	A good foundation for the optimization process is a well-defined low and high torque boundary region. 
	In a compression ignition engine the generated torque is roughly proportional to the injected fuel. 
	During the first run we measure the low torque region by limiting the \IF-interval to 6--10 $\text{mm}^3/\text{cycle}$, while the static actuators are set to values that support the creation 
	of low torque operation points; cf.\ Table~\ref{tab:actuators}.
	This run takes $6$ hours. 
	The low torque region is measured for a second time span of $6$ hours with slightly modified settings, i.e., \MT $= 10.0$ and \RP $= 405677 \operatorname{hPa}$. 
	We measure the low torque region twice because it represents, in particular, the situation 
	during startup, where low engine temperatures lead to unclean combustion, which causes high
	\HC and \CO emissions; cf.\ for example~\cite{Merkisz2015}.
	Hence, a detailed image of the engine behavior in this area of the operation field is desirable. 
	
	Afterwards, the \IF-interval is reset to 50--60 $\text{mm}^3/\text{cycle}$, so that the high torque regions can be measured.
	The static actuators are set to values that support high torque operation points; cf.\ Table~\ref{tab:actuators}.
	Again, the run takes $6$ hours.
	
	Finally, the lower and upper bound of \IF are removed, in order to get a picture of the remaining region of the solution map. 
	Therefore, the static actuators are set to midrange values; cf.\ Table~\ref{tab:actuators}. 
	This region is measured for $6$ hours, too.
	After $24$ hours, we get a coarse picture of the engine behavior.
	
	\begin{figure}
		\includegraphics[width=\textwidth]{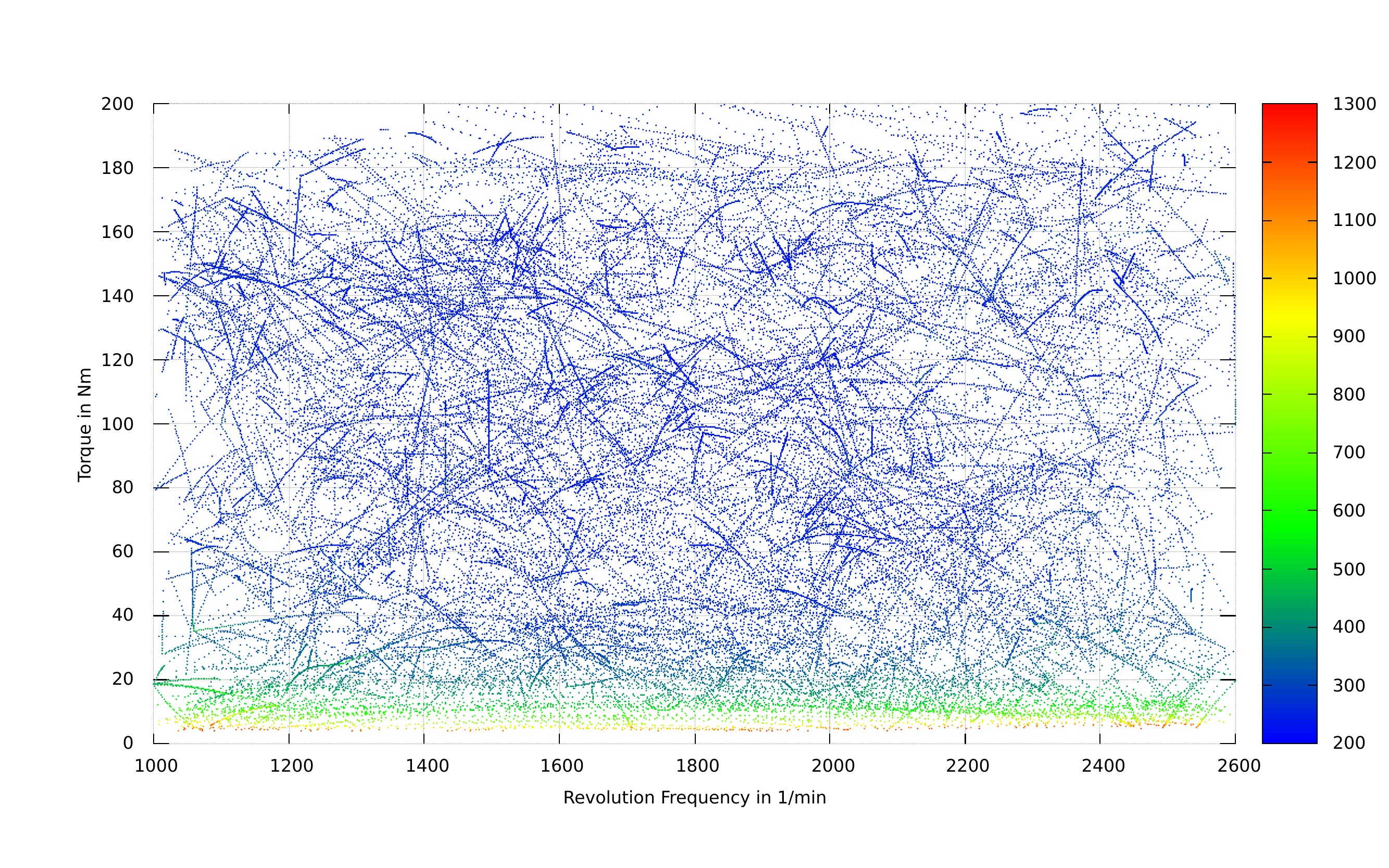}
		\caption{\label{fig:saturation}Operation field saturation after basic calibration, before closing of gaps.}
	\end{figure}
	
\begin{figure}[htbp]
		\begin{center}
			\begin{minipage}{0.9\textwidth}
				\centering
				\includegraphics[width=1.0\textwidth]{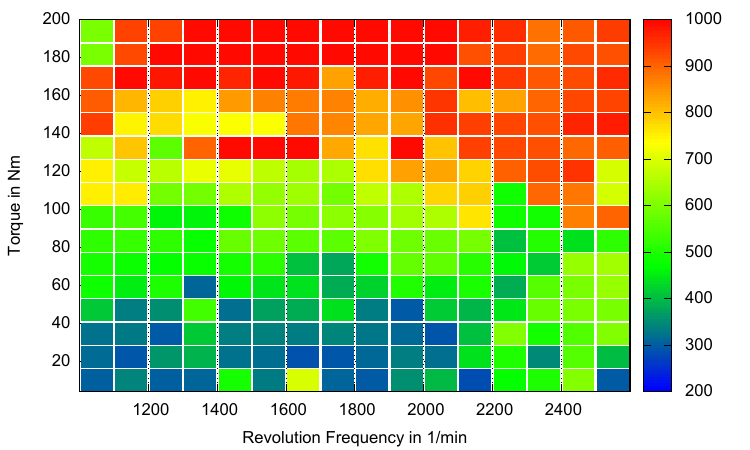}
				
			\end{minipage}
		\end{center}
		\vfill
		\begin{center}
			\begin{minipage}{0.9\textwidth}
				\centering
				\includegraphics[width=1.0\textwidth]{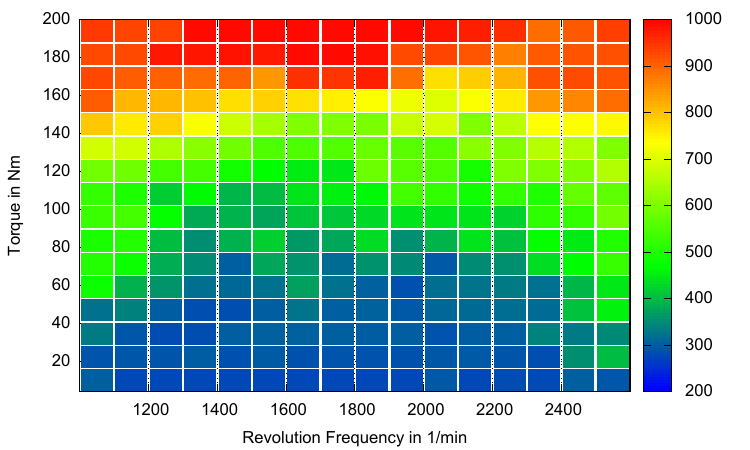}
			\end{minipage}
		\end{center}
		\caption{Top: Solution for actuator \AF without drivability constraint. Bottom: \AF settings of \euro 5 solution respecting the drivability constraint.}
		\label{fig:af-comp}
	\end{figure}

	As indicated above, the so-obtained data are used as a base set for a first run of 
	Algorithm~\ref{algo:calibration}, albeit with a deactivated emission constraint.
	Due to the omission of constraint~\eqref{eq:integral-emission-constraint}, the while loop
	exits immediately as a feasible solution $\SOLprelim$ to the ILP can always be found in this case. 
	Due to empty stacks and the drivability constraints, which are still active, there may still be gaps in
	$\SOLprelim$, though. 
	These are closed subsequently by interpoation-guided measurements; cf.\ Section~\ref{sec:hunt-mode}. 
	For the comparison of an engine map, i.e., a single component of a solution map, that respects the drivability constraint to one that does not, see Figure~\ref{fig:af-comp}.
	
	It took an additional $45.8$ hours to close all measurement holes. 
	Figure~\ref{fig:saturation} depicts the saturation of the operation field before this completion.
	The average fuel consumption for $100$ kilometers of the so-derived preliminary solution is $4.65$ liters in the \nedc.  
	Since this optimization does not take exhaust emissions into account, the cycle integrals of \CO, \HC and \NX are rather high: 3.68, 1.72 and 5.26, respectively, per \nedc. The solution map obtained in the first calibration phase serves as the basis of the full calibration. 
	
	\subsection{Full Calibration}\label{sec:full-calibration}
	
	During the full calibration phase the only static actuator is pilot injection. All others range over their whole respective domains. 
	With these presets, Algorithm~\ref{algo:calibration} runs automatically and without predefined
	time limit until it returns a solution map $\SOL$. 
	
	The optimization is performed, first, for the \nedc, then for the \rand cycle. An \nedc solution map that conforms to the \euro 4 norm is obtained after $68.51$ hours, one that conforms to \euro 5 after $88.41$ hours. 
	It turns out that the \nedc \euro 5 solution map conforms to all \euro 4 constraints with respect
	to the \rand cycle. 
	That is, it is also a \euro 4 solution map for the latter. 
	A \euro 5 solution map is obtained after $107.55$ hours of measurement.
	For the emission constraints corresponding to the different \euro norms, see Table~\ref{tab:euro-norms}.
	
	To lend some context to these numbers, we performed the same calibrations 
	using uniform grids. To derive a solution map of equivalent quality, i.e., a solution map conforming to all constraints given by the different \euro norms, the uniform grid approach consistently required 15 to 20 times as many measurements, whereby we again mean performed measurements,
	not stored ones. This translates into a real-world measurement time of weeks instead of days for
	a solution map of comparable quality. 
	
	\begin{remark}
		Since the AVL model only has a measurand for particle mass, but not for the number of particles, the determination of a \euro 6 solution map lies outside of the model's scope. However, the step from \euro 5 to \euro 6 merely adds an item to the list of emission constraints. This poses no fundamental challenge to our method which is scalable with respect to the number of pollutant limits. Naturally, adding further constraints will increase the measurement time though.
	\end{remark}
	
	\subsubsection*{Description of Figures}
	
\begin{figure}[htbp]
		\begin{center}
			\begin{minipage}{0.9\textwidth}
				\centering
				\includegraphics[width=1.0\textwidth]{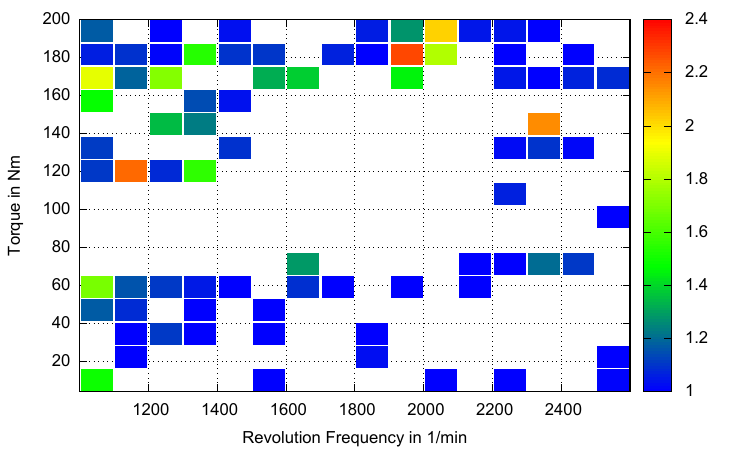}
			\end{minipage}
		\end{center}
		\caption{\footnotesize 
			Operation points where calibration for \rand cycle yields higher \NX emission than in \nedc. Colors represent the ratio of the emission values. }
		\label{fig:diff-nox}
	\end{figure}

	Figure~\ref{fig:diff-nox} displays the lower \NX output on the better part of the operation field for the \rand calibration in comparison to a calibration for the \nedc.
	
\begin{figure}[htbp]
		\begin{center}
			\begin{minipage}{0.9\textwidth}
				\centering
				\includegraphics[width=1.0\textwidth]{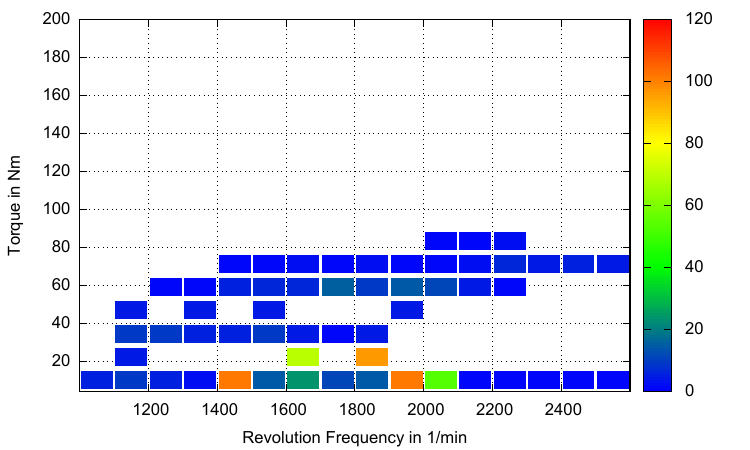}
			\end{minipage}
		\end{center}
		\vfill
		\begin{center}
			\begin{minipage}{0.9\textwidth}
				\centering
				\includegraphics[width=1.0\textwidth]{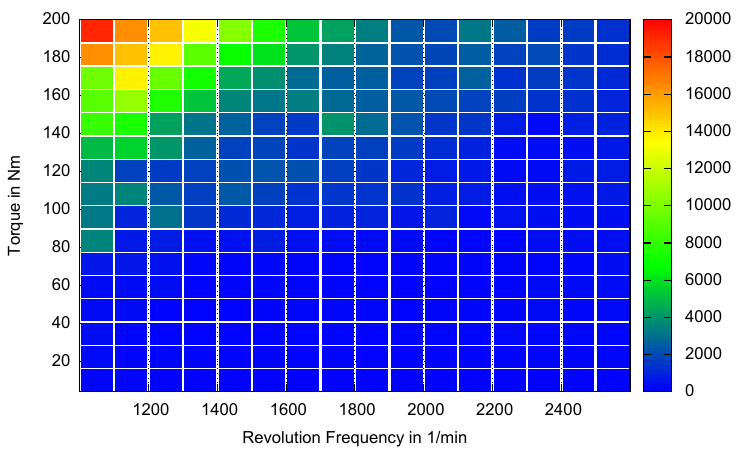}
				
			\end{minipage}
		\end{center}
		\vfill
		\begin{center}
			\begin{minipage}{0.9\textwidth}
				\centering
				\includegraphics[width=1.0\textwidth]{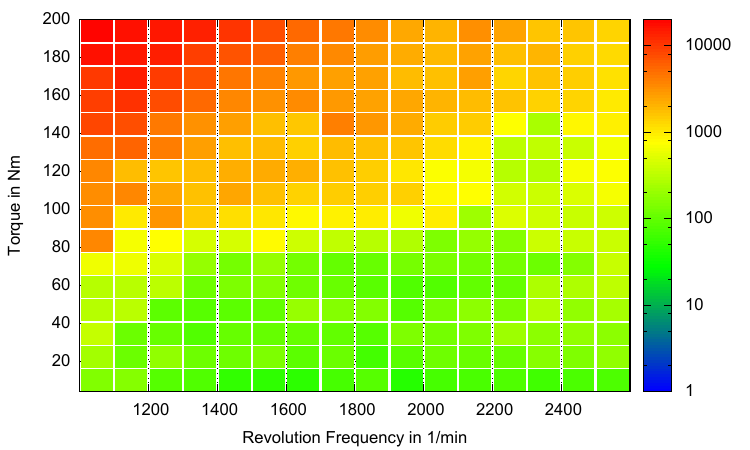}
				
			\end{minipage}
		\end{center}
		\caption{Top: Weighting table of the \nedc. Middle: Resulting \NX emissions (in $g/h$) after calibration for the \nedc. Bottom: \NX emissions with logarithmic scale.}
		\label{fig:nox-nedc}
	\end{figure}

\begin{figure}[htbp]
		\begin{center}
			\begin{minipage}{0.9\textwidth}
				\centering
				\includegraphics[width=1.0\textwidth]{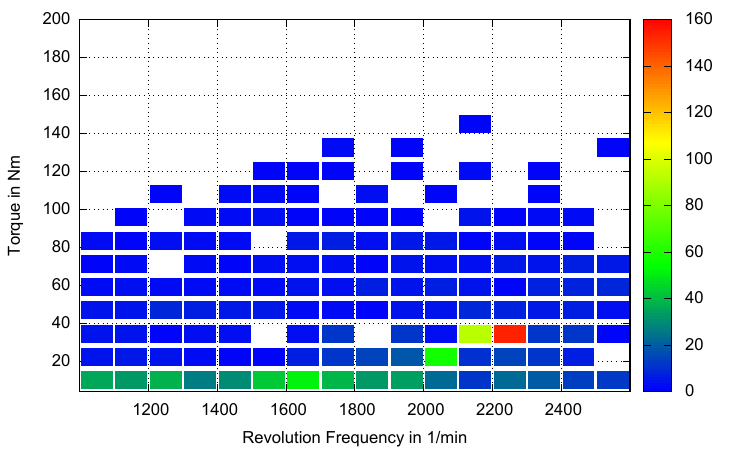}
			\end{minipage}
		\end{center}
		\vfill
		\begin{center}
			\begin{minipage}{0.9\textwidth}
				\centering
				\includegraphics[width=1.0\textwidth]{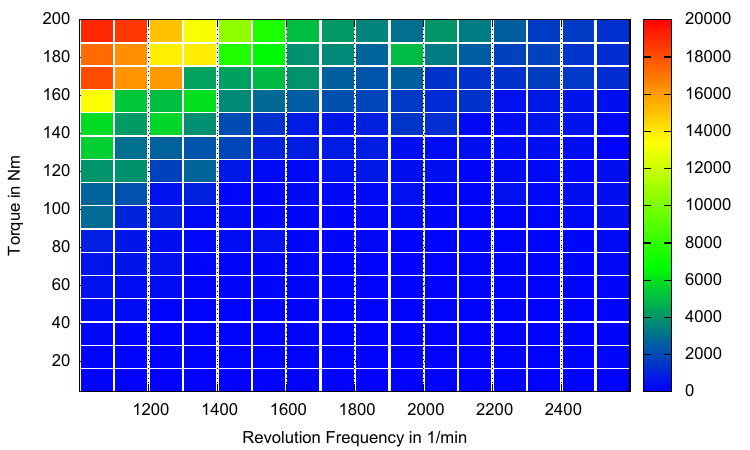}
				
			\end{minipage}
		\end{center}
		\vfill
		\begin{center}
			\begin{minipage}{0.9\textwidth}
				\centering
				\includegraphics[width=1.0\textwidth]{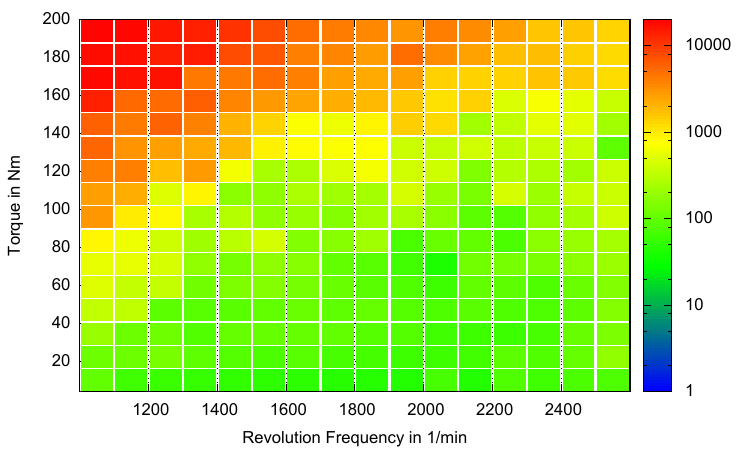}
				
			\end{minipage}
		\end{center}
		\caption{Top: Weighting table of the \rand cycle. Middle: Resulting \NX emissions (in $g/h$) after calibration for the \rand cycle. Bottom: \NX emissions on a logarithmic scale.}
		\label{fig:nox-random}
	\end{figure}

\begin{figure}[htbp]
		
		\begin{center}
			\begin{minipage}{0.9\textwidth}
				\centering
				\includegraphics[width=1.0\textwidth]{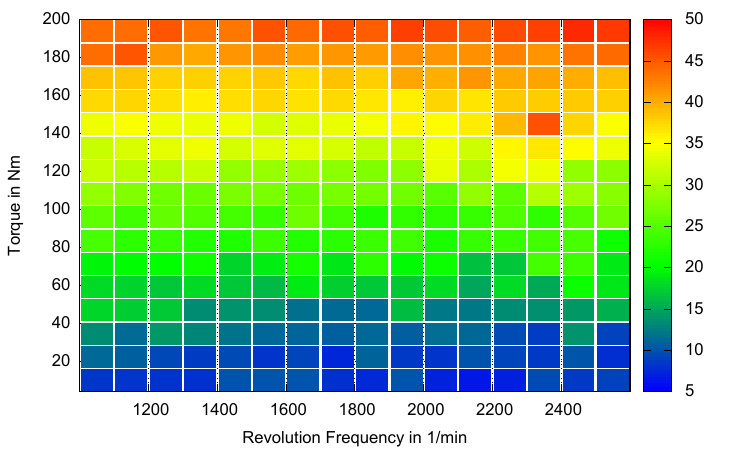}
				
			\end{minipage}
		\end{center}
		\vfill
		\begin{center}
			\begin{minipage}{0.9\textwidth}
				\centering
				\includegraphics[width=1.0\textwidth]{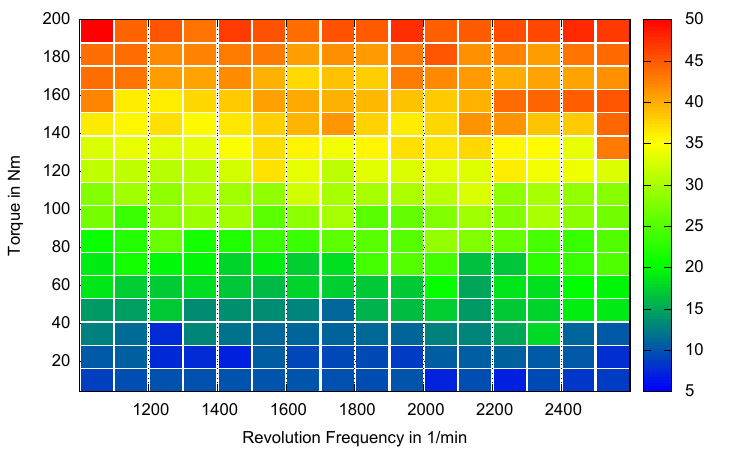}
				
			\end{minipage}
		\end{center}
		\vfill
		
		\begin{center}
			\begin{minipage}{0.9\textwidth}
				\centering
				\includegraphics[width=1.0\textwidth]{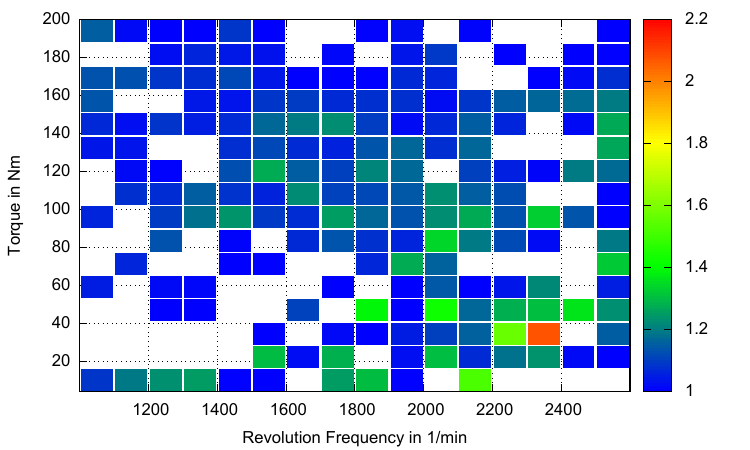}
				
			\end{minipage}
		\end{center}
		\caption{Actuator \IF after calibration for the \nedc (top) and \rand cycle (middle). Bottom: Operation points with higher values in calibration for \rand cycle than for \nedc. Colors represent ratio of consumption values.}
		\label{fig:diff-fl}
	\end{figure}

	Figures~\ref{fig:nox-nedc} and~\ref{fig:nox-random} illustrate a typical feature of engine calibrations fitted to specific driving cycles. These essentially subdivide the solution map into two parts. One covered by the cycle, where emissions are optimized, and one where they are not. This leads to emission values being 20--100 times higher on the part of the solution map that is not covered by the driving cycle. The \nedc covers only a fraction of the 
	operation field, while the \rand cycle covers more than half of it.
	In Figure~\ref{fig:diff-fl} one can observe that the solution map for the \rand cycle displays 10--15\,\% higher values for specific fuel consumption on most points of the operation field than the corresponding operation points of a solution for the \nedc.
	One can, of course, enforce the selection of data points with lower emission values on the part of the operation field not covered by a driving cycle as well, but at the expense of a higher fuel consumption.

	%
	%
	%
	%
	%
	%

	\section{Conclusion}
	
	In this article, we have described a semi-automatic approach to
	calibrate and optimize internal combustion engines with several
	actuators and sensors. The side constraints are limits given by
	safety or technical requirements, bounding the variation speed of actuators and
	ensuring emission bounds on given driving cycles. Our method automatically
	performs refinements of measurements, thus focusing the effort of the measurements around 
	regions of strongly nonlinear or even nonsmooth behavior of the engine. 
	That is, it automatically identifies neighborhoods 
	of anomalous engine behavior and maps them in appropriate detail.
	For the so-obtained data an optimal calibration solution is computed. 
	
	The output of the algorithm is a solution map $\SOL$ which consists of actual measurements
	and thus reflects the exact behavior of the engine for the given settings, as opposed
	to indirectly derived actuator setting obtained via modeling or interpolation. 
	This results in improved values for pollutant emission and fuel consumption near strongly nonlinear or even nonsmooth regions of the admissible domain.
	
	In our experiments, we demonstrated the practicability of the adaptive meshing methodology, 
	showing a significant speed-up in the measurement time (from weeks down to days) in comparison to uniform grids, without a loss of overall quality.
	Moreover, the resulting solution maps respect the emission constraints of \euro 4 and 5 norms. 
	We would like to stress that in our method it is easy to take into account further emission constraints such as, e.g., the number of emitted particles for the \euro 6 norms. 
	
	The next interesting step will be to test the method on an actual combustion engine. It is 
	our expectation that the experimental findings of this work will transfer well to the 
	real-life setting which is the engine test bench. 
	
	\section*{acknowledgements}
          We thank two anonymous reviewers for their comments that helped to improve the presentation of the paper.

	
	\bibliographystyle{spmpsci}
	\bibliography{engine}
	
\end{document}